\documentclass[final]{siamart171218}

\usepackage{amsfonts,MnSymbol,graphicx}
\usepackage[numbers]{natbib}
\usepackage{multirow,booktabs}

\usepackage{algorithmicx} 
\usepackage{algpseudocode}
\Crefname{ALC@unique}{Line}{Lines} 

\usepackage{subcaption} 

\usepackage{siunitx} 

\usepackage{accents}
\newcommand{\ubar}[1]{\underaccent{\bar}{#1}}

\DeclareMathOperator{\diag}{diag}

\newcommand*{\urho}{\ubar{\rho}} 
\newcommand*{\orho}{\bar{\rho}} 

\newcommand*{\unu}{\ubar{\nu}} 
\newcommand*{\onu}{\bar{\nu}} 

\newcommand*{\urr}{\ubar{r}} 
\newcommand*{\orr}{\bar{r}} 

\newcommand*{\ual}{\ubar{\alpha}} 
\newcommand*{\oal}{\bar{\alpha}} 

\newcommand*{\umu}{\ubar{\mu}} 
\newcommand*{\omu}{\bar{\mu}} 
\newcommand*{\uq}{\ubar{q}} 
\newcommand*{\oq}{\bar{q}} 
\newcommand*{\tg}{\tilde{g}} 
\newcommand*{\uh}{\ubar{h}} 
\newcommand*{\oh}{\bar{h}} 
\newcommand*{\rhop}{\rho'} 
\newcommand*{\rhopp}{\rho''} 
\newcommand*{\umup}{\umu'} 
\newcommand*{\umupp}{\umu''} 
\newcommand*{\omup}{\omu'} 
\newcommand*{\omupp}{\omu''} 

\newcommand{\RR}{\mathbb{R}}

\newcommand*{\tolsym}[1]{\varepsilon_{\scriptscriptstyle\rm #1}}
\newcommand*{\epsNWT}{\tolsym{NWT}}
\newcommand*{\epsPBM}{\tolsym{PBM}}
\newcommand*{\epsIP}{\tolsym{IP}}
\newcommand*{\epsDOC}{\tolsym{DOC}}

\newcommand*{\epsMINRES}{\tolsym{MR}}
\newcommand*{\res}{{\rm res}}
\newcommand*{\tres}{\widetilde{\res}}
\newcommand*{\resPBM}{\tres_{PBM}}
\newcommand*{\resIP}{\tres_{IP}}

\newcommand*{\altunderline}[1]{\mkern 2mu\underline{\mkern-0.5mu #1 \mkern-4.5mu}\mkern 6mu}
\newcommand*{\altoverline}[1]{\mkern 6mu\overline{\mkern-4.5mu #1 \mkern-0.5mu}\mkern 2mu}
\newcommand*{\UNu}{\altunderline{N}}
\newcommand*{\ONu}{\altoverline{N}}
\newcommand*{\URho}{\altunderline{P}}
\newcommand*{\ORho}{\altoverline{P}}

\date{}
\begin{document}
\bibliographystyle{plainnat}

\title{On barrier and modified barrier multigrid methods\\
for 3d topology optimization\thanks{This work has been partly supported by Fondation math\'{e}matique Jaques Hadamard FMJH/PGMO Project No2017-0088 "Multi-level Methods in Constrained Optimization".
}}
\author{Alexander Brune\thanks{School of Mathematics, University of Birmingham, Edgbaston,
	Birmingham B15 2TT, UK}
\and
Michal Ko\v{c}vara\thanks{School of Mathematics, University of
    Birmingham, Edgbaston, Birmingham B15 2TT, UK, and Institute of Information Theory
and Automation, Academy of Sciences of the Czech Republic, Pod
vod\'arenskou v\v{e}\v{z}\'{\i}~4, 18208 Praha 8, Czech Republic}}

\maketitle

\begin{abstract}
One of the challenges encountered in optimization of mechanical structures, in particular in what is known as topology optimization, is the size of the problems, which can easily involve millions of variables.
A basic example is the minimum compliance formulation of the variable thickness sheet (VTS) problem, which is equivalent to a convex problem. We propose to solve the VTS problem by the Penalty-Barrier Multiplier (PBM) method,  introduced by R.\ Polyak and later studied by Ben-Tal and Zibulevsky and others. The most computationally expensive part of the algorithm is the solution of linear systems arising from the Newton method used to minimize a generalized augmented Lagrangian. We use a special structure of the Hessian of this Lagrangian to reduce the size of the linear system and to convert it to a form suitable for a standard multigrid method. This converted system is solved approximately by a multigrid preconditioned MINRES method. The proposed PBM algorithm is compared with the optimality criteria (OC) method and an interior point (IP) method, both using a similar iterative solver setup. We apply all three methods to different loading scenarios. In our experiments, the PBM method clearly outperforms the other methods in terms of computation time required to achieve a certain degree of accuracy.
\end{abstract}



\section{Introduction}\label{sec:Intro}

The goal of topology optimization is to find an optimal geo\-metry of a solid body that maximizes its performance under certain boundary conditions, by determining an optimal distribution of material in a predefined design domain. It has many applications in industry, such as in mechanical and electrical engineering. The main challenge is the high computational cost of solving large-scale systems that arise from numerical methods to solve PDEs on high-resolution meshes. A basic example of topology optimization is the minimum compliance problem, where the deformation energy of an elastic body under prescribed loading and boundary conditions is to be minimized, given an amount of material. Relating the local stiffness of the body linearly to the continuous material distribution and employing a finite element discretization leads to the so-called \emph{variable thickness sheet} (VTS) problem
\begin{equation}
\label{eq:to_intro}
\begin{aligned}
  &\min_{\rho\in\RR^m\!,\,u\in\RR^n} \frac{1}{2}f^\top u\\
  &\mbox{subject to}\\
  &\qquad K(\rho) u = f\\
  &\qquad \sum_{i=1}^m \rho_i = V\\
  &\qquad \rho_i\geq \urho_i, \quad i=1,\ldots,m \\
  &\qquad \rho_i\leq \orho_i, \quad i=1,\ldots,m  \; ,
\end{aligned}
\end{equation}
where $K(\rho) = \sum_{i=1}^m \rho_i K_i$, with $K_i\in\RR^{n\times n}$, is the stiffness matrix and $f\in \RR^n$ is the load vector of the finite element equilibrium equations.
The design variable $\rho$ is commonly referred to as the \emph{density}, while the vector $u$ represents the nodal displacements.
We assume that $K_i$ are symmetric and positive semidefinite and that
$\sum_{i=1}^m K_i$ is sparse and positive definite. We also assume that the
volume $V\in\RR$ and the lower and upper bounds $\urho\in\RR^m_+$ and $\orho\in\RR^m_+$ are chosen such that the problem is strictly feasible. This implies $\orho>\urho$, among other things. While problem \eqref{eq:to_intro} is not itself convex, it is equivalent to a convex problem; see \cite{ben1996hidden} and Theorem~\ref{th:equiv} below. For a more detailed derivation of the VTS problem and a comprehensive treatment of the theory and applications of topology optimization, see for example \cite{bendsoe-sigmund}.

The minimum compliance problem has been studied extensively. Still, it is the subject of ongoing research as higher design detail calls for higher mesh resolution, which in turn makes the problem more computationally demanding. \citeauthor{Aage_2017}, for example, performed topology optimization on a model with more than one billion elements \cite{Aage_2017}. The bottleneck of algorithms for topology optimization is usually the solution of large linear systems. Direct solvers are not a viable option, due to their computational complexity and demand on computer memory, and iterative, most typically Krylov type solvers, are given preference. Since their convergence behavior highly depends on the condition number of the system matrix, preconditioning plays a vital role. The multigrid method, introduced by \citeauthor{Brandt_1977} as a solver for boundary-value problems \cite{Brandt_1977}, has become popular as a means to precondition the system by employing it inside the iterative solvers. As early as \citeyear{maar-schulz}, \citet{maar-schulz} proposed a conjugate gradient (CG) method preconditioned by multigrid for topology optimization. Similar solvers were used in \cite{amir} and \cite{MK_Mohammed_2016}. In \cite{Aage_2017}, the authors chose a multi-layered algorithm involving two types of Krylov solvers and the geometric as well as algebraic multigrid method. We refer the reader to \cite{briggs2000multigrid} for a comprehensive introduction to the multigrid method.

Beyond the issue of efficiently solving the linear systems arising within each iteration of the optimization algorithm, the total number of such iterations required to reach the optimal solution---and thus the choice of optimization method---also affects the overall time-efficiency of the algorithm. The most commonly used methods for the minimum compliance problem are the \emph{optimality criteria} (OC) method, see \cite{bendsoe-sigmund}, and the \emph{method of moving asymptotes} (MMA) \cite{svanberg}. One of the advantages of the OC method is that it is relatively simple to implement; see in particular \cite{top88}. To our knowledge, global convergence results exist only for the MMA and it is often the algorithm of choice in commercial software or large-scale applications, such as that described in \cite{Aage_2017}. Both of these methods, however, usually rely on heuristics for their stopping criteria and, in practice, display a very similar rate of convergence.

A possible alternative to the aforementioned methods is the \emph{interior point} (IP) method. It has become increasingly popular in the past twenty to thirty years, particularly for convex optimization \cite{Wright_1997:IP}. Its theoretical advantage over the OC method or the MMA for convex problems lies in its rate of convergence, especially for convex quadratic problems such as \eqref{eq:to_intro}. \citet{maar-schulz} used an IP algorithm for 2D topology optimization. In \cite{jarre-kocvara-zowe}, \citeauthor{jarre-kocvara-zowe} proposed an IP method for truss topology optimization. This was later extended in \cite{MK_Mohammed_2016} to large 2D VTS problems, where it outperformed the OC method, in terms of both iterations and overall CPU time required to achieve optimality to within a certain precision. In one part of our paper, we build on this work and further improve the algorithm to apply it to large-scale 3D problems. The approach is described in Section~\ref{sec:IP} and results of some examples are presented in Section~\ref{sec:Num}. 

Going from 2D to 3D is by no means straightforward. The largest examples in \cite{MK_Mohammed_2016} were based on nine regular refinements of a very coarse, e.g. $2\times 2$, mesh. This resulted in 262\,144 finite elements and 526\,338 degrees of freedom (components of the displacement vector $u$). Such a problem could still be solved on a standard laptop. If we used the same refinement level in a 3D example starting with a $2\times 2\times 2$ coarse mesh, we would end up with a problem with more than 134 million finite elements and 405 million degrees of freedom. Moreover, while the stiffness matrix in 2D typically has 18 non-zero elements per row, in 3D problems this number typically goes up to 81 non-zeros, i.e., the stiffness matrix is considerably denser. All this makes much greater demands on the numerical linear algebra used in the optimization algorithm.

A common problem with IP methods is the ill-conditioning of the system as the iterates approach the optimal solution. This leads to an increase in solver iterations which can make the algorithm nonviable. A class of methods that aims to counteract this problem while otherwise following a strategy similar to that of the IP method, is the class of \emph{penalty-barrier multiplier} (PBM) methods. They were first introduced in \cite{ben1997penalty}, building on the modified barrier methods proposed by \citeauthor{polyak1992modified} in \cite{polyak1992modified}. As part of the larger class of augmented Lagrangian methods, they have one particular convergence property which sets them apart from IP methods. The latter involve a sequence of barrier parameters which needs to tend to 0 for convergence to the optimal solution, this being the cause of the increasing ill-conditioning; the former feature a penalty parameter for which there exists a value larger than 0 such that the method still converges to the optimal solution. See, for example, \cite[Corollary 6.15]{Stingl_2006} for a result specific to penalty-barrier methods. PBM methods have been successfully applied to convex problems and semidefinite problems in topology optimization \cite{pennon-iter}. In Section~\ref{sec:MGNR} of this paper, a penalty-barrier method for \eqref{eq:to_intro} is introduced. In contrast to the IP method, the PBM method does not stay in the strict interior of the feasible region. This poses a problem with regard to the positive definiteness of $K(\rho)$, which depends on $\rho_i$ being strictly positive for all $i=1,\dots,m$. We circumvent this problem by applying the PBM method to the dual of \eqref{eq:to_intro}. The theoretical background for this is covered in Section~\ref{sec:PrimalDualVTS}. The PBM approach described in Section~\ref{sec:MGNR} is applied to several examples in Section~\ref{sec:Num}, in order to compare it to the IP method from Section~\ref{sec:IP}, as well as to the OC method, which is briefly described in Section~\ref{sec:OC}.

Lastly, a remark on notation: throughout this paper, we use $e_i$ to denote the $i$-th canonical unit vector and $e$ to denote the vector $(1,\dots,1)^\top$ of appropriate dimension.

\section{Dual VTS problem}\label{sec:PrimalDualVTS}
Consider the variable thickness sheet problem \eqref{eq:to_intro}. Following \cite{Ben-Tal1993,kocvara2017truss} in the context of equivalent formulations for truss topology optimization, we can formulate a dual problem to \eqref{eq:to_intro}:
\begin{equation}
\label{eq:dual}
\begin{aligned}
&\min_{u\in\RR^n\!,\,\alpha\in\RR,\,\unu,\,\onu\in\RR^m}
\alpha V-f^\top u -\urho^\top \unu+\orho^\top \,{\onu}\\
&\mbox{subject to}\\
&\qquad \frac{1}{2}u^\top K_iu\leq \alpha-\unu_i+{\onu}_i, \quad i=1,\ldots,m  \\
&\qquad {\unu}_i\geq 0,\quad i=1,\ldots,m\\
&\qquad {\onu}_i\geq 0,\quad i=1,\ldots,m\,.
\end{aligned}
\end{equation}

\begin{theorem}\label{th:equiv}
	Problems \eqref{eq:to_intro} and \eqref{eq:dual} are equivalent in
	the following sense:
	\begin{itemize}
		\item[(i)] If one problem has a solution then also the other problem
		has a solution and
		$$
		\min\eqref{eq:to_intro} = \min\eqref{eq:dual}\,.
		$$
		\item[(ii)] Let
		$(u^*,\alpha^*,\unu^*,\onu^*)$ be a solution
		to \eqref{eq:dual}. Further, let $\tau^*$ be the vector of Lagrangian
		multipliers for the inequality constraints associated with this
		solution. Then $(u^*,\tau^*)$ is a solution of \eqref{eq:to_intro}.
		Moreover,
		$$
		\unu^*_i\onu^*_i = 0,\quad i=1,\ldots,m.
		$$
		\item[(iii)] Let $(u^*,\rho^*)$ be a solution of \eqref{eq:to_intro}.
		Further, let $\urr^*$ and $\orr^*$ be the Lagrangian
		multipliers associated with the lower and upper bounds on $\rho$,
		respectively, and let $\lambda^*$ be the multiplier for the volume constraint. Then $(u^*,\lambda^*,\urr^*,\orr^*)$
		is a solution of \eqref{eq:dual}.
	\end{itemize}
\end{theorem}
\begin{proof}
	We will first write
	\eqref{eq:to_intro} equivalently as
	\begin{equation}
	\min_{\stackrel{\urho_i \leq \rho_i \leq
			\orho_i}{\sum_{i=1}^m \rho_i = V}} \max_{u\in\RR^n}\ f^\top u -
	\frac{1}{2} u^\top K(\rho)u \,.\label{compliance_b_var}
	\end{equation}
	Indeed, as $K(\rho)$ is by assumption positive semidefinite, the necessary and sufficient optimality condition for the inner maximization problem is $K(\rho)u=f$ and, using this, the optimal value of the maximization problem is 
	$\frac{1}{2}f^\top u$.
	Problem \eqref{compliance_b_var} is convex (actually linear) and bounded in $\rho$ and
	concave in $u$, so we can switch ``max'' and ``min'' (see, e.g.,
	\cite{ekeland-temam}) to get an equivalent problem:
	$$
	\max_{u\in\RR^n}\inf_{\stackrel{\urho_i \leq \rho_i \leq
			\orho_i}{\sum_{i=1}^m \rho_i = V}}  f^\top u - \frac{1}{2} u^\top K(\rho)u
	\,.
	$$
	Due to our assumption of strict feasibility, there exists a Slater
	point for the feasible set of the inner (convex) optimization
	problem, so we may replace it
	by its Lagrangian dual. The Lagrangian multipliers for the inequalities
	will be denoted by $\urr\in\RR^m_+$
	and $\orr\in\RR^m_+$, that for the volume equality constraint by $\lambda\in\RR$:
	\begin{equation} \label{eq:qcqp_th1}
	\max_{u\in\RR^n}\max_{\stackrel{\scriptstyle
			\lambda\in\RR}{\urr\in\RR^m_+,
			\orr\in\RR^m_+}}\inf_{\rho\in\RR^m_+}\ 
	f^\top u - \frac{1}{2} u^\top K(\rho)u + \lambda(\sum_{i=1}^m \rho_i - V)
	-\urr^\top (\rho-\urho)
	+\orr^\top (\rho-\orho)  \,.
	\end{equation}
    We can include the non-negativity constraint on $\rho$ in the inner-most optimization problem because we know that the solution to \eqref{eq:qcqp_th1} satisifies $\rho\ge\urho\ge 0$.
	
	Now regard the dual problem \eqref{eq:dual}. It can equivalently be formulated as the following min-max problem, 
	using a partial Lagrangian function with multiplier $\tau\in\RR^m$:
	$$
	\min_{\stackrel{\scriptstyle u\in\RR^n}{\stackrel{\scriptstyle\alpha\in\RR}
			{\unu\in\RR^m_+,\onu\in\RR^m_+}}}
	\max_{\tau\in\RR^m_+}\ 
	\alpha V - f^\top u -\unu^\top \urho + \onu^\top \orho
	+ \sum_{i=1}^m \tau_i(\frac{1}{2} u^\top K_iu - \alpha+\unu_i - \onu_i)
	$$
	which can be rearranged further to give
	\begin{equation} \label{eq:qcqp_th2}
	\min_{\stackrel{\scriptstyle u\in\RR^n}{\stackrel{\scriptstyle\alpha\in\RR}
			{\unu\in\RR^m_+,\onu\in\RR^m_+}}}
	\max_{\tau\in\RR^m_+}\ 
	\frac{1}{2} u^\top K(\tau)u - f^\top u + \alpha(V-\sum_{i=1}^m \tau_i)
	+\unu^\top (\tau-\urho)
	-\onu^\top (\tau-\orho)  \,.
	\end{equation}
	Identifying $\tau$, $\alpha$, $\unu$, and $\onu$ with $\rho$, $\lambda$, $\urr$, and $\orr$, respectively, and changing the sign of the objective
	function (and thus changing ``max'' to ``min'' and ``min'' to
	``max''), we can
	see that \eqref{eq:qcqp_th1} and \eqref{eq:qcqp_th2} are equivalent.
	For later reference, note that the multiplier $\tau$ of the dual problem corresponds to the primal variable $\rho$, the density.
	
	The second part of (ii) is obvious from the fact that
	$\unu$ and $\onu$ are multipliers for the lower
	and upper bounds, so only one of them can be positive (only one
	bound can be active) for each component.
\end{proof}

Notice that \eqref{eq:dual} is a convex optimization problem, as $K_i$ are positive semidefinite.

We finish this section with another formulation of the dual VTS problem that allows us to easily compute the duality gap (this formulation was first derived in \cite{Ben-Tal1993}).
\begin{theorem}
	Problem \eqref{eq:dual} is equivalent to an
	unconstrained nonsmooth problem
	\begin{align} \label{eq:qcqp_c}
	&\max_{u\in\RR^n,\alpha\in\RR} -\alpha V+f^\top u + \sum_{i=1}^m
	\min\{(\alpha-\frac{1}{2} u^\top K_iu)\urho_i,(\alpha-\frac{1}{2} u^\top
	K_iu)\orho_i\}
	\end{align}
	in the following sense:
	\begin{itemize}
		\item[(i)] $\min\eqref{eq:dual}= - \max\eqref{eq:qcqp_c}$;
		\item[(ii)] Let
		$(u^*,\alpha^*,\unu^*,\onu^*)$ be a solution
		of \eqref{eq:dual}. Then $(u^*,\alpha^*)$ is a solution of~\eqref{eq:qcqp_c}. 
		Conversely, every solution $(u^*,\alpha^*)$
		of \eqref{eq:qcqp_c} is a part of~a solution of \eqref{eq:dual}.
	\end{itemize}
\end{theorem}
\begin{proof}
	We will show that \cref{eq:dual} and \cref{eq:qcqp_c} are equivalent reformulations of each other. Introducing an auxiliary variable $s\in\RR^m$, problem
	\eqref{eq:qcqp_c} can be directly re-written as
	\begin{align*}
	&\max_{\scriptstyle u\in\RR^n,\alpha\in\RR,s\in\RR^m}
	-\alpha V+f^\top u + \sum_{i=1}^m s_i \\
	&\mbox{\rm subject to}\\
	&\qquad (\alpha-\frac{1}{2}u^\top K_iu)\orho_i\geq s_i,
	\quad i=1,\ldots,m \\
	&\qquad (\alpha-\frac{1}{2}u^\top K_iu)\urho_i\geq s_i,
	\quad i=1,\ldots,m  \,.
	\end{align*}
	The constraints in the above problem can be written as
	$$
	(\alpha-\frac{1}{2}u^\top K_iu)\geq
	\max\{\frac{s_i}{\urho_i},\frac{s_i}{\orho_i}\}\,,
	\quad i=1,\ldots,m \,.
	$$
	Noting that $\orho>\urho\geq 0$, we define
	\begin{alignat*}{2}
	& \unu_i = \frac{s_i}{\urho_i}\,,\
	\onu_i=0\,, && \quad\text{if}\quad
	\frac{s_i}{\urho_i}>\frac{s_i}{\orho_i} > 0\\
	& \unu_i = 0\,,\ \onu_i=-\frac{s_i}{\orho_i}\,, 
	&& \quad\text{if}\quad
	\frac{s_i}{\urho_i}\leq\frac{s_i}{\orho_i} \leq 0\,.
	\end{alignat*}
	Then the above set of constraints can also be written as
	$$
	(\alpha-\frac{1}{2}u^TK_iu)\geq
	\unu_i - \onu_i  \quad i=1,\ldots,m \,.
	$$
	Obviously, these $\unu_i,\onu_i$ also satisfy the non-negativity constraints. Lastly, we can reformulate the objective function to match \eqref{eq:dual}, since
	$$
	\sum_{i=1}^m  {\urho_i}\unu_i - \sum_{i=1}^m
	{\orho_i}\onu_i =
	\sum_{i:\frac{s_i}{\urho_i}>\frac{s_i}{\orho_i}}
	{\urho_i}\frac{s_i}{\urho_i} + \sum_{i:
		\frac{s_i}{\urho_i}\leq\frac{s_i}{\orho_i}}
	{\orho_i}\frac{s_i}{\orho_i} = \sum_{i=1}^m s_i\,.
	$$
	We switch the sign of the objective function and claims (i) and (ii) follow.
\end{proof}

Assume that $(u,\alpha)$ is a feasible point in the dual problem \eqref{eq:dual} such that there exist $\rho$ satisfying $K(\rho)u= f$ and $(\rho,u)$ is feasible in the primal problem \eqref{eq:to_intro}.
We then have the following formula for the duality gap:
\begin{equation}
\label{eq:gap}
\begin{aligned}
\delta(u,\alpha) :=& \; \min{\eqref{eq:to_intro}} - \max{ \eqref{eq:qcqp_c} }  \\
=& -\frac{1}{2}f^\top u + 
\alpha V - \sum_{i=1}^m \min\left\{\urho_i(\alpha-\frac{1}{2}u^\top K_iu),\orho_i(\alpha-\frac{1}{2}u^\top K_iu)\right\}\,.
\end{aligned}
\end{equation}

\section{The penalty-barrier multilplier method for topology optimization}\label{sec:MGNR}
In this section, we describe the class of Penalty-Barrier Multiplier (PBM) algorithms and their application to the VTS problem. This class of algorithms was originally developed and analyzed by R.~Polyak under the name Modified Barrier algorithms; see, among others, \cite{polyak1988smooth,polyak1992modified}. These methods are defined for a class of ``modified" barrier functions; a particular choice of a function leads to a particular algorithm. Ben-Tal and Zibulevsky \cite{ben1997penalty}  analyzed one such choice that proved to be computationally very efficient; see also \cite{kocvara2003pennon}. The PBM method was first applied to topology optimization problems in \cite{kocvara1998mechanical}.

\subsection{Penalty-barrier multiplier methods}
Consider a generic convex constraint optimization problem 
$$
  \min_x \{\mathfrak{f}(x)\mid \mathfrak{g}_i(x)\leq 0,\  i=1,\ldots,m\}\,.
$$ 
The idea of NR is to replace the inequalities by scaled inequalities $p_i\varphi\left(\frac{\mathfrak{g}_i(x)}{p_i}\right)\leq 0$ with a penalty function $\varphi$ and a penalty parameter $p_i>0$. 
Here, $\varphi$ is a strictly increasing, twice differentiable, real-valued, strictly convex function 
with dom $\varphi = (- \infty, b), \ 0 < b \leq \infty$,
which has the following properties:
\begin{itemize}
	\item[$(\varphi 1)$] $\qquad \varphi (0) = 0$ 
	\item[$(\varphi 2)$] $\qquad \varphi^{\prime} (0) = 1$ 
	\item[$(\varphi 3)$] $\qquad{\displaystyle \lim_{s \rightarrow b}} \ 
	\varphi^{\prime} (s) = \infty $ 
	\item[$(\varphi 4)$] $\qquad 
	{\displaystyle \lim_{s \rightarrow - \infty}} 
	\varphi^{\prime} (s) = 0 $. 
\end{itemize}
\smallskip\noindent
Then the ``penalized" problem 
\begin{equation}\label{eq:penalized}
  \min_x \{\mathfrak{f}(x)\mid p_i\varphi\left(\frac{\mathfrak{g}_i(x)}{p_i}\right)\leq 0,\ i=1,\ldots,m\}
\end{equation}
remains convex and has the same feasible set and thus the same solution as the original one. We formulate a standard Lagrangian function of the penalized problem that can be considered an augmented Lagrangian function of the original problem:
\begin{equation}\label{eq:lagr}
  {\cal L}(x,\mu; p) = \mathfrak{f}(x) + \sum_{i=1}^m \mu_i p_i \varphi\left(\frac{\mathfrak{g}_i(x)}{p_i}\right)\,.
\end{equation}

At each iteration of the NR, we minimize the augmented Lagrangian with respect to $x$
\begin{align}
\mbox{\it Step 1.} &\qquad x^{k+1} \approx \arg\min_x {\cal L}(x, \mu^k; p^k) 
\label{eq:105}\\
\intertext{and update the multipliers and the penalty parameter:}
\mbox{\it Step 2.} &\qquad \mu_i^{k+1} = \mu_i^k \varphi^{\prime} \left(\frac{\mathfrak{g}_i(x^{k+1})} {p_i^k}\right) 
\label{eq:106} \\
\mbox{\it Step 3.} &\qquad p_i^{k+1} = \pi p_i^k \ .  
\label{eq:107} 
\end{align}
Here $\pi<1$ is a penalty updating factor. The meaning of the ``$\approx$" sign in Step~1 is that the unconstrained minimization problem is only solved approximately, until $\|\nabla_{\!x}\,{\cal L}(x,\mu;p)\|\leq\varepsilon$, where $\varepsilon$ is some prescribed tolerance.

For more details on the NR methods, analysis and numerical performance, see the references above.

In Step 1 we need to solve, approximately, an unconstrained optimization problem. For this, we will use the Newton method. Therefore,
we will need formulas for the gradient and Hessian of ${\cal L}$ with respect to the primal variable $x$:
\begin{equation}\label{eq:lagr_grad}
  \nabla_{\!x}\,{\cal L}(x,\mu;p) = \nabla_{\!x}\, \mathfrak{f}(x) + \sum_{i=1}^m \mu_i \varphi' \left(\frac{\mathfrak{g}_i(x)}{p_i}\right)\nabla_{\!x}\, \mathfrak{g}_i(x)
\end{equation}
and
\begin{equation}
\label{eq:lagr_hess}
\begin{aligned}
\nabla^2_{\!xx}\,{\cal L}(x,\mu;p) = &\nabla^2_{\!xx}\, \mathfrak{f}(x) + \sum_{i=1}^m \frac{\mu_i}{p_i} \varphi'' \left(\frac{\mathfrak{g}_i(x)}{p_i}\right)\nabla_{\!x}\, \mathfrak{g}_i(x)(\nabla_{\!x}\, \mathfrak{g}_i(x))^\top \\
&+ \sum_{i=1}^m \mu_i \varphi' \left(\frac{\mathfrak{g}_i(x)}{p_i}\right) \nabla^2_{\!xx}\,\mathfrak{g}_i(x)\,.
\end{aligned}
\end{equation}
Note that, due to the convexity of the penalized problem \eqref{eq:penalized}, the Hessian of ${\cal L}$ is positive semidefinite for any arguments $x\in\RR^n$, $\mu\in\RR^m_+$.

\citeauthor{ben1997penalty} \cite{ben1997penalty} analyzed one particular choice of the penalty function $\varphi$ defined as follows:
\begin{equation}
\varphi_{\hat{\tau}} (\tau) = \left \{ 
\begin{array}{ll} 
	\tau + \frac{1}{2} \, \tau^2 & \tau \geq \hat{\tau} \\[0.2cm] 
	- (1+ \hat{\tau})^2 \log \left ( \frac{1+ 2 \hat{\tau} -\tau}
	{1 + \hat{\tau}} \right)
	+ \hat{\tau} + \frac{1}{2} \hat{\tau}^2 \quad & \tau < \hat{\tau}  \ . 
\end{array} \right .  
\label{103.1}
\end{equation}
By setting $\hat{\tau} = - \frac{1}{2}$,
we get a pure (not shifted) logarithmic branch.
As this function combines properties of the quadratic penalty function and the logarithmic barrier function, it is called a penalty-barrier function and the resulting algorithm a penalty-barrier multiplier method. This method proved to be very efficient and we will use it to solve the dual VTS problem.

\subsection{PBM for the dual VTS problem}
Let us now apply the PBM method to the dual problem \eqref{eq:dual}. The augmented Lagrangian for this problem is defined as
\begin{equation}
\label{eq:aug_lagr}
\begin{aligned}
  {\cal L}(u,\alpha,\unu, \onu, \rho,\umu, \omu) = \ &
 \alpha V-f^\top u -\urho^\top\unu+\orho^\top\onu \\
  &+ \sum_{i=1}^m \rho_i p_i \varphi\left(\frac{1}{p_i}(\frac{1}{2}u^\top K_i u - \alpha + \unu_i  - \onu_i)\right) \\
  &+ \sum_{i=1}^m \umu_i \uq_i \varphi\left(\frac{-\unu_i}{\uq_i}\right) 
  + \sum_{i=1}^m \omu_i \oq_i \varphi\left(\frac{-\onu_i}{\oq_i}\right)
\end{aligned}
\end{equation}
with Lagrangian multipliers $\rho\in\RR^m$, $\umu\in\RR^m$  and $\omu\in\RR^m$ and penalty parameters $p\in\RR^m$, $\uq\in\RR^m$ and $\oq\in\RR^m$. 

To simplify the notation, let us define the aggregate variable 
$$
\xi := (u,\alpha,\unu,\onu)
$$
and vectors of penalized constraints as
\begin{align*}
 \tg_i(\xi)&=\tg_i(u,\alpha,\unu,\onu) := \varphi\left(\frac{1}{p_i}(\frac{1}{2}u^\top K_i u - \alpha +\unu_i - \onu_i)\right), \quad i=1,\ldots,m \,,\\
\uh_i(\xi)&=\uh_i(u,\alpha,\unu,\onu) := \varphi\left(\frac{-\unu_i}{\uq_i}\right), \quad i=1,\ldots,m\,,\\
\oh_i(\xi)&=\oh_i(u,\alpha,\unu,\onu) := \varphi\left(\frac{-\onu_i}{\oq_i}\right), \quad i=1,\ldots,m\,.
\end{align*}
Let $s_i(\xi)$ denote the argument of $\varphi(\cdot)$ in the definition of $\tg_i(\xi)$ above. In the following, the notation $\tg'_i(\xi)$ will be understood as $\varphi^\prime(s_i(\xi))$, rather than a composite derivative of $\varphi(s_i(\xi))$ with respect to $\xi$. We define $\uh'_i(\cdot)$ and $\oh'_i(\cdot)$ analogously, as well as $\tg''(\cdot)$, $\uh''(\cdot)$ and $\oh''(\cdot)$.

According to \eqref{eq:lagr_grad}, the gradient of the augmented Lagrangian with respect to the aggregate variable $\xi$  is
\begin{equation}
\label{eq:lagr_grad_a}
\begin{aligned}
&\nabla_{\!\xi}\,{\cal L}(\cdot) 
= \begin{bmatrix}\res_1\\\res_2\\\res_3\\\res_4\end{bmatrix}
=\begin{bmatrix} -f\\V\\\urho\\\orho\end{bmatrix} +
\sum_{i=1}^m \rho_i\tilde{g}'_i(\xi)\begin{bmatrix} K_iu\\-1\\e_i\\-e_i\end{bmatrix} +
\sum_{i=1}^m \umu_i\uh'_i(\xi)\begin{bmatrix} 0\\0\\-e_i\\0\end{bmatrix} +
\sum_{i=1}^m \omu_i\oh'_i(\xi)\begin{bmatrix} 0\\0\\0\\-e_i\end{bmatrix}\,.
\end{aligned}
\end{equation}
To further simplify the notation, we define
\begin{align*}
\rhop_i = \rhop_i(\xi) &:= \rho_i \tg'_i(\xi) \,, &
\umup_i = \umup_i(\xi) & := \umu_i \uh'_i(\xi) \,, &
\omup_i = \omup_i & := \omu_i \oh'_i(\xi) \,, \\[0.2em]
\rhopp_i = \rhopp_i(\xi) &:= \dfrac{\rho_i}{p_i}\tg''_i(\xi) \,, &
\umupp_i = \umupp_i & := \dfrac{\umu_i}{\uq_i} \uh''_i(\xi) \,, &
\omupp_i = \omupp_i & := \dfrac{\omu_i}{\oq_i} \oh''_i(\xi) \,. \\
\end{align*}

By \eqref{eq:lagr_hess}, the Hessian of the augmented Lagrangian will take the form
\begin{equation}\label{eq:lagr_hess_a}
\nabla^2_{(\!u,\alpha,\unu,\onu)^2}\,{\cal L}(\cdot) = 
\begin{bmatrix}
H_{11}&H_{12}&H_{13}&H_{14}\\
H_{12}^\top &H_{22}&H_{23}&H_{24}\\
H_{13}^\top &H_{23}^\top &H_{33}&H_{34}\\
H_{14}^\top &H_{24}^\top &H_{34}^\top&H_{44}
\end{bmatrix} 
\end{equation}
where
\begin{align*}
  H_{11} &= 
  \displaystyle\sum_{i=1}^m \rhopp_i K_iu u^\top K_i^\top  
  + \displaystyle\sum_{i=1}^m \rhop_i K_i,\quad H_{11}\in\RR^{n\times n}\\
  H_{12}&=-\displaystyle\sum_{i=1}^m \rhopp_i K_iu ,\quad H_{12}\in\RR^{n\times 1}\\
  H_{13}& = \left[\rhopp_1 K_1u,\ldots,\rhopp_m K_mu\right],
  \quad H_{13}\in\RR^{n\times m}\\
  H_{14}& = \left[-\rhopp_1 K_1u,\ldots,-\rhopp_m K_mu\right],
  \quad H_{14}\in\RR^{n\times m}\\
  H_{22}&=\displaystyle\sum_{i=1}^m\rhopp_i,\quad H_{22}\in\RR\\
  H_{23}&=\left[ -\rhopp_1,\ldots, -\rhopp_m\right],\quad H_{23}\in\RR^{1\times m}\\
  H_{24}&=\left[\rhopp_1,\ldots,\rhopp_m\right],\quad H_{24}\in\RR^{1\times m}\\
  H_{33}&=\diag(\rhopp_1+\umupp_1,\ldots,\rhopp_m+\umupp_m),\quad H_{33}\in\RR^{m\times m}\\
  H_{34}&=\diag(-\rhopp_1,\ldots,-\rhopp_m),\quad H_{34}\in\RR^{m\times m}\\
  H_{44}&=\diag(\rhopp_1+\omupp_1,\ldots,\rhopp_m+\omupp_m),\quad H_{44}\in\RR^{m\times m}\,.
\end{align*}
By \eqref{eq:106}, the Lagrange multipliers in the PBM algorithm are never equal to zero. Hence, the matrices $H_{33}, H_{34}, H_{44}$ are diagonal and positive or negative definite, so we can easily calculate the following inverse of the lower right block of the Lagrangian, which is in turn a block diagonal matrix:
$$
\begin{bmatrix}H_{33} &H_{34}\\H_{34}^\top &H_{44}\end{bmatrix}^{-1} =
\begin{bmatrix}H_{33}^{-1}+H_{33}^{-1}H_{34}ZH_{34}^\top H_{33}^{-1} &-H_{33}^{-1}H_{34}Z\\
-ZH_{34}^\top H_{33}^{-1} &Z\end{bmatrix}
$$
with $Z=(H_{44}-H_{34}^\top H_{33}^{-1}H_{34})^{-1}$. We will require this inverse further below.

Observe that the matrix $H_{11}$ has the same sparsity structure as the ``unscaled" stiffness matrix $\displaystyle\sum_{i=1}^m K_i$. Indeed, the only non-zero components of the vector $K_iu$ are those corresponding to indices of non-zero elements of $K_i$, hence $K_i$ has the same sparsity structure as $(K_iu)(K_iu)^\top $. For this reason, the matrices $H_{13} H_{13}^\top $ and $H_{14} H_{14}^\top $ have the same sparsity structure as $H_{11}$ and thus $\displaystyle\sum_{i=1}^m K_i$. This property extends to any matrices $H_{13} D H_{13}^\top $ and $H_{14} D H_{14}^\top $, where $D$ is a diagonal matrix. 

\begin{figure}[h]
		\begin{center}
			\begin{subfigure}{0.45\textwidth}
                \includegraphics[width=\hsize]{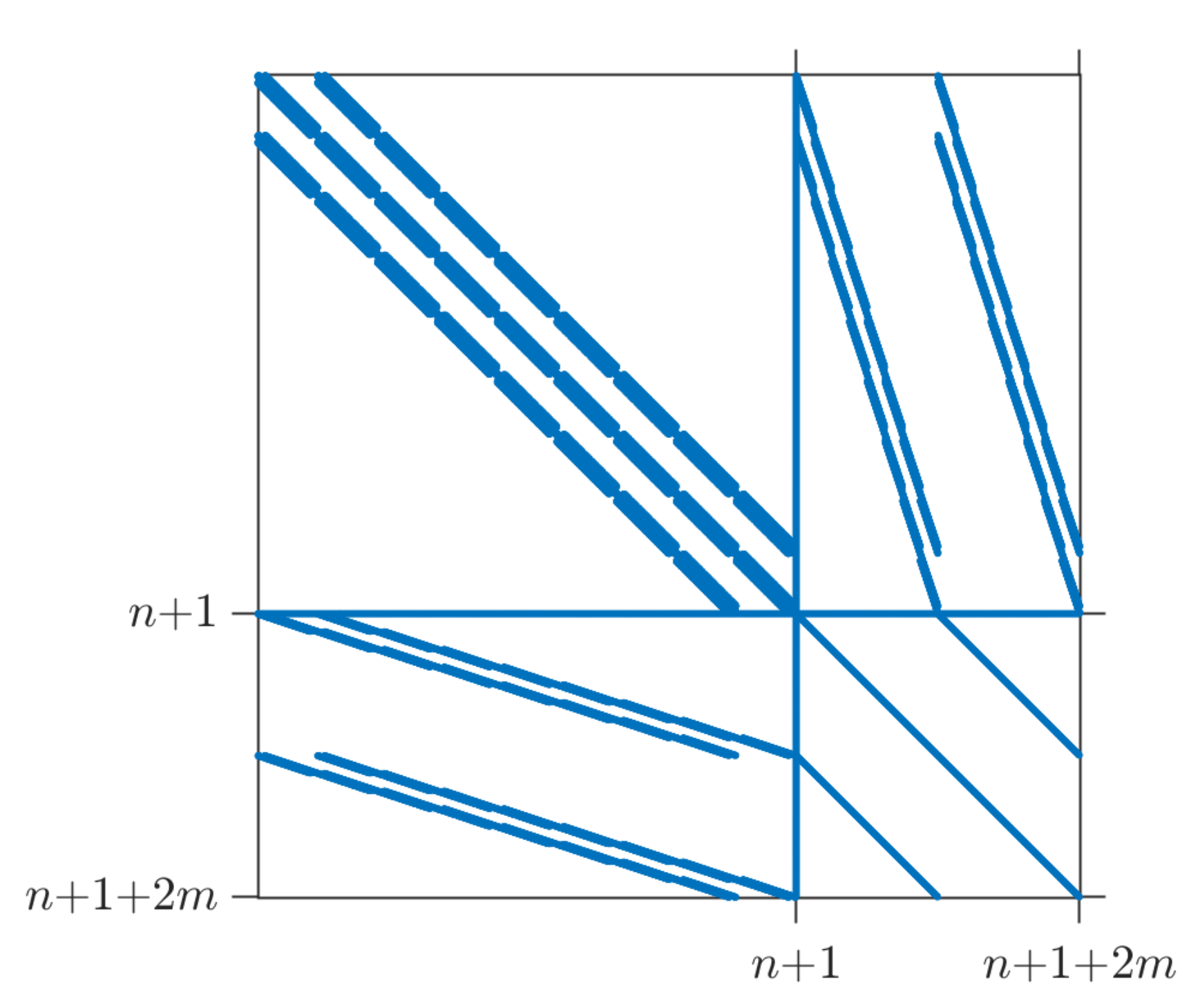}

				\caption{$\nabla^2 \cal L$}
			\end{subfigure} \hspace{2em}
			\begin{subfigure}{0.45\textwidth}
                \includegraphics[width=\hsize]{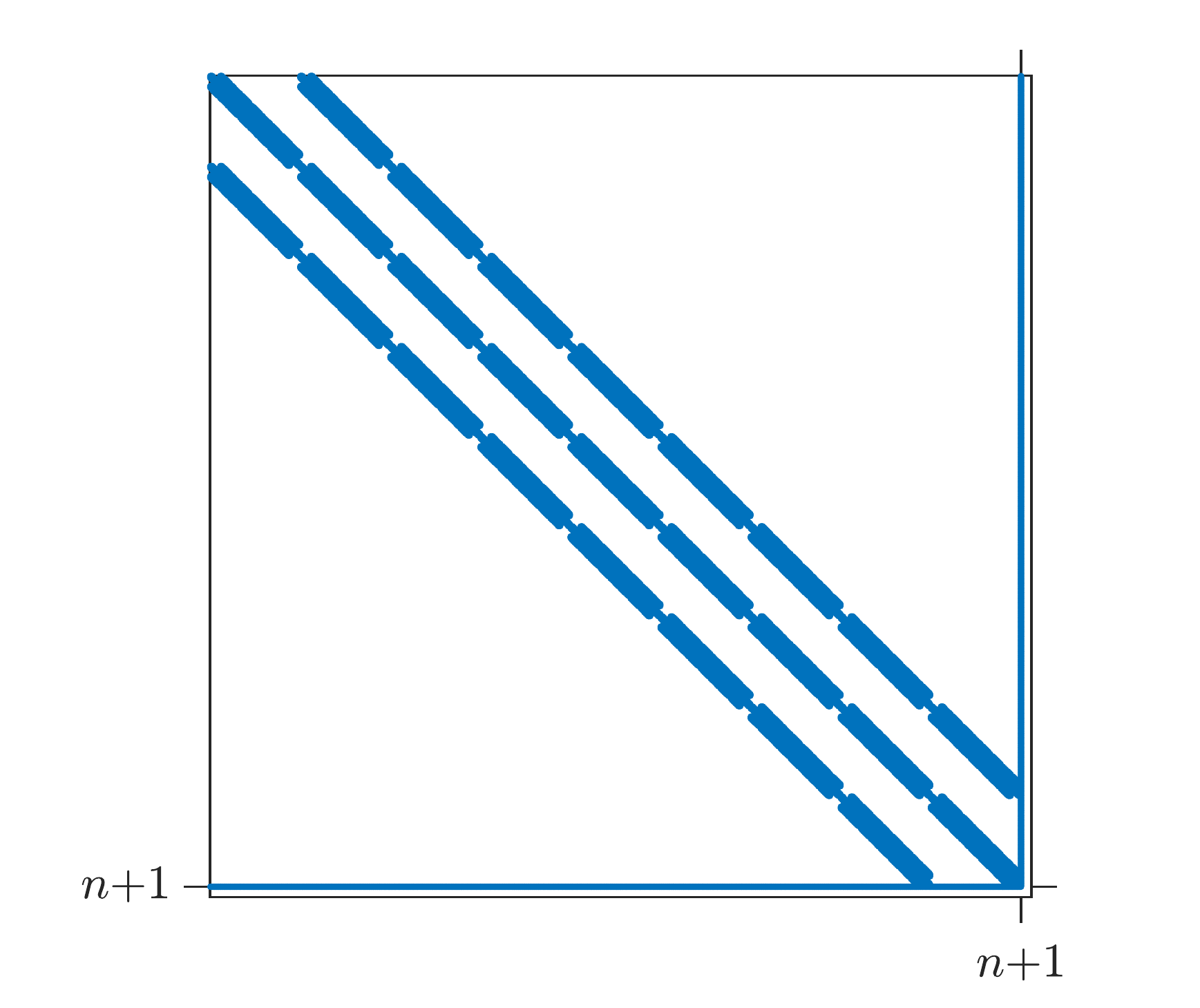}
				\caption{$S$}
				\label{subfig:1right}
			\end{subfigure}
		\end{center}	
	\caption{\label{fig:1}Typical sparsity structure of the Hessian of the augmented Lagrangian for the dual topology optimization problem (left) and its Schur complement (right).}
\end{figure}

We now calculate the following Schur complement matrix
\begin{equation}\label{eq:S}
  S = \begin{bmatrix} H_{11}&H_{12}\\ H_{12}^\top &H_{22} \end{bmatrix} 
  - \begin{bmatrix}H_{13}& H_{14}\\H_{23} &H_{24}\end{bmatrix}
  \begin{bmatrix}H_{33} &H_{34}\\H_{34}^\top &H_{44}\end{bmatrix}^{-1}
  \begin{bmatrix}H_{13}^\top &H_{23}^\top\\ H_{14}^\top&H_{24}^\top \end{bmatrix}
  \; \in\RR^{(n+1)\times (n+1)} \, .
\end{equation}
By the previous considerations, the principal $n\times n$ submatrix of $S$ has the same sparsity structure as the stiffness matrix $\displaystyle\sum_{i=1}^m K_i$; the last row and column of $S$ are full. Figure~\ref{fig:1} shows typical examples of the sparsity structure of the Hessian of the augmented Lagrangian $\nabla^2_{\!\xi\xi}\,{\cal L}(\cdot)$ in \eqref{eq:lagr_hess_a} and the Schur complement matrix $S$.

The first step of the PBM algorithm is to solve approximately the unconstrained minimization problem
$$
  \min_{u,\alpha,\unu,\onu} {\cal L}(u,\alpha,\unu, \onu, \rho,\umu,\omu)
$$
by the Newton method. In every step of the Newton method, we have to solve the system of linear equations
$$
\nabla^2_{(\!u,\alpha,\unu,\onu)^2}\,{\cal L} (u,\alpha,\unu,\onu, \rho,\umu,\omu)\cdot (\Delta u,\Delta\alpha,\Delta\nu)
= -\nabla_{(\!u,\alpha,\unu,\onu)}\,{\cal L}(u,\alpha,\unu,\onu, \rho,\umu,\omu)\,,
$$
where $(\Delta u,\Delta\alpha,\Delta\nu)$ is the Newton increment and $\Delta \nu:=(\Delta\unu,\Delta\onu)$.
Equivalently, according to the above development, we can instead solve the reduced system
\begin{equation}\label{eq:system}
  S \; \begin{bmatrix}
  \Delta u\\\Delta\alpha
  \end{bmatrix}= \mbox{\it rhs} \,,
\end{equation}
where, by \eqref{eq:lagr_grad_a},
\begin{align*}
  \mbox{\it rhs} = 
  &-\begin{bmatrix} -f\\V\end{bmatrix} -
  \sum_{i=1}^m \rhop_i\begin{bmatrix} K_iu\\-1\end{bmatrix}\\[0.2em]
  &+ \begin{bmatrix}H_{13}&H_{14}\\ H_{23}&H_{24}\end{bmatrix} \begin{bmatrix}H_{33} &H_{34}\\H_{34}^\top &H_{44}\end{bmatrix}^{-1}  \left( \begin{bmatrix}\urho\\ \orho\end{bmatrix}
  + \sum_{i=1}^m \rhop_i \begin{bmatrix} e_i \\ -e_i \end{bmatrix} 
  + \sum_{i=1}^m \umup_i \begin{bmatrix} e_i \\ 0 \end{bmatrix}
  + \sum_{i=1}^m \omup_i \begin{bmatrix} 0 \\ -e_i \end{bmatrix} \right) \,.
\end{align*}
Recall that the dual problem \eqref{eq:dual} is convex, hence the Hessian of ${\cal L}$ is positive semidefinite and, consequently, so is the Schur complement $S$. 

The remaining component $\Delta \nu$ can be reconstructed from the solution to \eqref{eq:system} as follows:
\begin{equation}
\label{eq:deltanu}
\begin{aligned}
  \Delta \nu = \;  - \begin{bmatrix}H_{33} &H_{34}\\H_{34}^\top &H_{44}\end{bmatrix}^{-1}
 & \left( \; \begin{bmatrix}\urho\\\orho\end{bmatrix} 
  + \sum_{i=1}^m \rhop_i \begin{bmatrix}e_i\\-e_i\end{bmatrix} 
  + \sum_{i=1}^m \umup_i \begin{bmatrix}e_i\\0\end{bmatrix}  
  + \sum_{i=1}^m \omup_i \begin{bmatrix}0\\-e_i\end{bmatrix} \right. \\[0.2em]
  & + \left. \begin{bmatrix}H_{13}^\top \\ H_{14}^\top\end{bmatrix} \Delta u 
  + \begin{bmatrix}H_{23}^\top \\ H_{24}^\top\end{bmatrix} \Delta\alpha\right)\,.
\end{aligned}
\end{equation}

After the augmented Lagrangian has been minimized, we check for convergence. For this, we use the duality gap $\delta(u,\alpha)$ in \eqref{eq:gap}, scaled by the dual objective function, henceforth denoted by $d(u,\alpha,\unu, \onu)$, as a measure of optimality. If convergence has not yet been achieved, the multipliers are updated, imposing the safeguard rule used in \cite{ben1997penalty}, followed by the penalty parameters. 

The PBM method is summarized in Algorithm~\ref{alg:pbm}. It employs the Newton method with backtracking line search using the Armijo rule; see Algorithm~\ref{alg:pbm_newton}. The stopping criterion for the Newton method uses the weighted residual term 
\begin{equation} \label{eq:resIP}
\resPBM = \dfrac{ \| \res_1 \|_2 }{ \| f \|_2 } + \dfrac{ | \res_2 | }{ V }
+ \dfrac{ \| \res_3 \|_2 }{ \| \urho \|_2 + \| \orho \|_2 }
\end{equation}
as a measure of feasibility. The stopping parameter is adjusted adaptively in each PBM iteration and this warrants some clarification. Setting the Newton method tolerance too low in early stages of the PBM method leads to an increase in Newton iterations and thus in computational time without significantly changing the convergence behavior of the PBM. A ``soft'' tolerance of 100 times the current optimality measure has proven to be a good choice. At the same time, however, we want to guarantee that the final solution has a certain degree of feasibility, which requires the system \eqref{eq:system} to be solved to a certain accuracy. For this reason, after the final PBM iteration, we run the Newton method one more time with decreased tolerance and then update $\rho$. For the sake of completeness, it should be noted that the solution $(u,\alpha,\unu,\onu)$ obtained by this additional call to the Newton method is not guaranteed to still satisfy the stopping criterion on Line \ref{algline:pbm_stopcrit} of Algorithm \ref{alg:pbm}. It is possible that it was previously only satisfied due to the inaccuracy of the solution\footnote{Note that $\delta(u,\alpha)$ is only a valid duality gap for \emph{feasible} solutions $u$ and $\alpha$.}. In the vast majority of our numerical experiments, however, this was not an issue and $|\delta(u,\alpha)/d(u,\alpha,\unu, \onu)|$ remained below the stopping parameter $\epsPBM$.

\begin{algorithm}
	\caption{PBM}	
	\label{alg:pbm}
	Let $0<\beta< 1$, $0<\gamma< 1$, $p_{\min},\uq_{\min},\oq_{\min}>0$, $\epsPBM>0$,
	$\epsNWT>0$ and $\epsNWT^{\min}>0$ be given.
	Choose initial vectors $(u,\alpha,\unu, \onu)$ and $(\rho,\umu,\omu)$.
	Set $p=\uq=\oq=e\in\RR^m$. \\[-\baselineskip]
	\begin{algorithmic}[1]
	\Repeat
		\State{Minimize the augmented Lagrangian \eqref{eq:aug_lagr} with respect to $(u,\alpha,\unu, \onu)$ by Algorithm~\ref{alg:pbm_newton} with stopping tolerance $\epsNWT$ \label{algline:pbm_newton}}
        \State{Compute the duality gap $\delta(u,\alpha)$ by \eqref{eq:gap} and the dual objective function value $d(u,\alpha,\unu, \onu)$}
        \If{$|\delta(u,\alpha)/d(u,\alpha,\unu, \onu)| < \epsPBM$} \label{algline:pbm_stopcrit}
        	\State STOP
        \EndIf
        \State Update the multipliers 
			\begin{align*}
			  \rho_i^+  &=  \rho_i \varphi^{\prime}\left(\frac{1}{p_i}(\frac{1}{2}u^\top K_i u - \alpha + \unu_i  - \onu_i)\right),\quad i=1,\ldots,m \\[0.5em]
			  \umu_i^+  &=  \umu_i \varphi^{\prime}\left(\frac{\unu_i}{\uq_i}\right), \quad
			  \omu_i^+  =  \omu_i \varphi^{\prime}\left(\frac{-\onu_i}{\oq_i}\right),\quad i=1,\ldots,m
		    \end{align*}
		\State If necessary, correct the multipliers such that
    		$$
    		  \beta\rho_i \leq \rho_i^+ \leq\frac{1}{\beta}\rho_i,\quad
    		  \beta\umu_i \leq \umu_i^+ \leq\frac{1}{\beta}\umu_i,\quad
    		  \beta\omu_i \leq \omu_i^+ \leq\frac{1}{\beta}\omu_i,\quad i=1,\ldots,m
    	    $$
    	    and set\ \ 
    	    $ \rho = \rho^+,\ \umu=\umu^+,\ \omu=\omu^+.$
    	    \label{algline:pbm_update_mult}
		\State{Update the penalty parameters
		    $$
	           p_i = \max\{\gamma\, p_i,p_{\rm min}\},\quad \uq_i = \max\{\gamma\, \uq_i,\uq_{\rm min}\},
	           \quad \oq_i = \max\{\gamma\, \oq_i,\oq_{\rm min}\},
	        $$
	        for $i=1,\ldots,m$
	        \State Update the stopping tolerance for Algorithm~\ref{alg:pbm_newton} 
	        $$
	          \epsNWT = \max \left\{\, \min \left\{ 100 \cdot \left|\dfrac{\delta(u,\alpha)}{d(u,\alpha,\unu, \onu)} \right| \,, \epsNWT \right\} \; ,
	          \; \epsNWT^{\min} \, \right\}
	        $$
        }
        \Until{convergence}
        \State Set $\epsNWT = 10\cdot\epsPBM$ and repeat Line \ref{algline:pbm_newton}
        \State Update $\rho$ as done in Line \ref{algline:pbm_update_mult}
	\end{algorithmic}
\end{algorithm}
Our choice of parameters in Algorithm~\ref{alg:pbm} was $\beta = \gamma = 0.3$, $p_{\min}=\uq_{\min}=\oq_{\min} = 10^{-8}$, $\epsNWT=1$ and $\epsNWT^{\rm min}=10^{-3}$. The initial values were $u=0$, $\alpha=1$, $\unu=\onu=e$, $\rho_i=V/m$, for all $i=1,\dots,m$, and $\umu=\omu=e$.

Note that Algorithm~\ref{alg:pbm_newton} is an inexact Newton method, which uses a preconditioned Krylov subspace method, as described later in Section~\ref{sec:MG}. Let us reiterate that the principal $n\times n$ submatrix of $S$ in \eqref{eq:system} has the same sparsity structure as the stiffness matrix $K(\rho)$. This will allow us in Section~\ref{sec:MG} to develop a  multigrid preconditioner using the standard prolongation/restriction operators for the stiffness matrix.
\begin{algorithm}
\caption{PBM NEWTON}	
\label{alg:pbm_newton}
Let vectors $(u,\alpha,\nu)$ and $(\rho,\umu,\omu)$ be given, using $\nu=(\unu,\onu)$.
Let $\epsNWT$ be given.
\begin{algorithmic}[1]
\Repeat
\State{Compute matrix $S$ from \eqref{eq:S} and the corresponding right hand side}
\State{Solve (approximately) the linear system \eqref{eq:system} to find $\Delta u, \Delta\alpha$}
\State{Compute $\Delta\nu$ from \eqref{eq:deltanu} with data $\Delta u, \Delta\alpha$}
\State{Perform backtracking line search with Armijo rule to find step length $\kappa$}
\State{Update $u,\alpha,\nu$: 
	$
	u = u + \kappa\Delta u,\quad
	\alpha = \alpha + \kappa\Delta \alpha,\quad
	\nu = \nu + \kappa\Delta \nu
	$
}
\If{$\resPBM < \epsNWT$}
\State {STOP}
\EndIf
\Until convergence
\end{algorithmic}
\end{algorithm}

\section{An Interior Point method for topology optimization}\label{sec:IP}
In this section, we describe the primal-dual IP method used to solve \eqref{eq:to_intro}. This involves deriving the linear system to be solved in each iteration and taking Schur complements of this system in order to obtain a system that, firstly, is symmetric positive definite and, secondly, displays a structure that allows a straightforward application of the multigrid method as a preconditioner. In this, we follow \cite{MK_Mohammed_2016}. Many features of the algorithm proposed in that reference had to be changed to make it more performant and viable for 3D problems. Therefore, we include all details of the algorithm. We do not recapitulate the basics of primal-dual IP methods and instead refer the reader to \cite{Wright_1997:IP}, to name just one standard piece of literature.

Some notation from the previous section will be reused below for variables that serve a similar purpose. However, the primal and dual variables $\rho$, $u$, $\alpha$, $\unu$ and $\onu$ have the same meaning in both sections. This is worthwhile to note because it means that the results from the PBM method described in the previous section and the IP method described below are directly comparable.

\subsection{Primal-dual Interior Point method for the VTS problem}
We start by setting up the KKT conditions for the VTS problem \eqref{eq:to_intro}. Note that the problem exhibits a ``hidden convexity'', i.e., it is not itself a convex problem but is equivalent to a different, convex problem \cite{ben1996hidden}. The strict feasibility, given for \eqref{eq:to_intro} by design---see Section \ref{sec:Intro}---translates to this equivalent problem. Hence, the Slater condition is satisfied and the KKT conditions are necessary and sufficient optimality conditions. They are given by the constraint equations in \eqref{eq:to_intro} and the equations below.
\begin{align*}
\dfrac{1}{2} u^\top K_i u + \alpha + \unu_i - \onu_i &= 0 \,, \quad i=1,\dots,m \\
(\rho_i - \urho_i) \unu_i &= 0 \,, \quad i=1,\dots,m \\
(\orho_i - \rho_i) \onu_i &= 0 \,, \quad i=1,\dots,m \; .
\end{align*}
Note that in the above, the Lagrange multipliers for the equilibrium equation constraint $K(\rho)u=f$ have already been eliminated, taking advantage of the fact that the minimum compliance problem is self-adjoint. This means that, due to our choice of objective function, the aforementioned multipliers also satisfy the equilibrium equation---with the the right-hand side only differing by a constant factor. Hence, we can directly identify them with $u$. See, for example, \cite{bendsoe-sigmund} for details.

The complementarity conditions for the lower and upper bound constraints, i.e., the second and third lines in the system above, are now perturbed by replacing $0$ by barrier parameters $r>0$ and $s>0$, respectively.
The resulting system of equations needs to be solved for fixed $r,s$ in each iteration of the IP algorithm. This is done approximately by performing one iteration of the Newton method. We get the following residual function for the Newton method:
\begin{align*}
\res(u,\alpha,\rho,\unu,\onu) \; = \; \begin{bmatrix} \res_1 \\ \res_2 \\ \res_3 \\ \res_4 \\ \res_5 \end{bmatrix} \; =
\; & \begin{bmatrix} -f \\ -V \\ 0 \\ -r \, e \\ -s \, e \end{bmatrix} + 
\sum_{i=1}^m \rho_i \begin{bmatrix} K_i u \\ 1 \\ 0 \\ 0 \\ 0  \end{bmatrix} + 
\sum_{i=1}^m \dfrac{1}{2}u^\top K_i u \begin{bmatrix} 0 \\ 0 \\ e_i \\ 0 \\ 0 \end{bmatrix} \\[0.5em]
+ & \sum_{i=1}^m \alpha \begin{bmatrix} 0 \\ 0 \\ e_i \\ 0 \\ 0 \end{bmatrix} +
\sum_{i=1}^m \unu_i \begin{bmatrix} 0 \\ 0 \\ e_i \\ (\rho_i-\urho_i) e_i \\ 0 \end{bmatrix} + 
\sum_{i=1}^m \onu_i \begin{bmatrix} 0 \\ 0 \\ e_i \\ (\orho_i-\rho_i) e_i \\ 0 \end{bmatrix}
\end{align*}
Next, we obtain the derivative of the residual function as the block matrix
\begin{equation}
\label{eq:gradres}
\nabla_{(u,\alpha,\rho,\unu,\onu)} \res(\cdot) = 
\begin{bmatrix}
    K(\rho) & 0 & B(u) & 0 & 0 \\
    0 & 0 & e^\top & 0 & 0 \\
    B(u)^\top & e & 0 & I & -I \\
    0 & 0 & \UNu & \URho & 0 \\
    0 & 0 & -\ONu & 0 & \ORho
\end{bmatrix} \; ,
\end{equation}
where $I\in\RR^{m\times m}$ is the identity matrix and we use the notation
{
\setlength{\jot}{0.8em}
\begin{align*}
B(u) &= \left[ K_1 u, \dots, K_m u \right] \,,\\
\UNu &= \diag( \unu ) \,,\qquad \ONu = \diag( \onu ) \,,\\
\URho &= \diag( \rho - \urho ) \,, \qquad
\ORho = \diag( \orho - \rho ) \,. 
\end{align*}%
}%
The system matrix $\nabla \res$ in \eqref{eq:gradres} is indefinite. Similar to the procedure in Section \ref{sec:MGNR}, we can reduce the above system to a positive definite one. We do this in two steps. First, we construct the Schur complement of $\nabla \res$ with respect to its invertible lower right block $\begin{bmatrix} \URho & 0 \\ 0 & \ORho \end{bmatrix}$. We then in turn form the Schur complement of the result with respect to its lower right block; see \cite{MK_Mohammed_2016} for details. This leaves us with the matrix
\begin{equation}
\label{eq:S_IP}
S =
\begin{bmatrix} K(\rho) & 0 \\ 0 & 0 \end{bmatrix} +
\begin{bmatrix} B(u) \\ e^\top \end{bmatrix}
\left( \URho^{-1} \UNu + \ORho^{-1} \ONu \right) ^{-1}
\begin{bmatrix} B(u)^\top & e \end{bmatrix}
\;\, \in\RR^{(n+1)\times(n+1)} \,.
\end{equation}
This matrix is positive definite as long as $\rho$ is strictly feasible and $\unu,\onu > 0$. Recall that $(K_i u)(K_i u)^\top$ has the same sparsity structure as $K_i$. Hence, the matrix $S$ in \eqref{eq:S_IP} has the same sparsity structure as that in \eqref{eq:S} in the previous section.

In each iteration of the IP method, we approximately solve the nonlinear system
\[
\res(u,\alpha,\rho,\unu,\onu) = 0
\]
by performing one iteration of Newton's method. Instead of solving the Newton system
\[
\nabla_{(u,\alpha,\rho,\unu,\onu)} \res(u,\alpha,\rho,\unu,\onu) \cdot (\Delta u, \Delta \alpha, \Delta \rho, \Delta \unu, \Delta \onu) = 
- \res (u,\alpha,\rho,\unu,\onu) \, ,
\]
we solve the equivalent system
\begin{equation}
\label{eq:IP_system}
S \; \begin{bmatrix} \Delta u \\ \Delta \alpha \end{bmatrix} = rhs \, ,
\end{equation}
where, according to the above reduction of the system,
\begin{align*}
rhs = & -\begin{bmatrix} -f \\ -V \end{bmatrix} - \sum_{i=1}^m \rho_i \begin{bmatrix} K_i u \\ 1 \end{bmatrix} \\[0.5em]
& - \begin{bmatrix} B(u) \\ e^\top \end{bmatrix} \left( \URho^{-1} \UNu + \ORho^{-1} \ONu \right) ^{-1}
\left(
\res_3 + \URho^{-1} \res_4 - \ORho^{-1} \res_5
\right) \, .
\end{align*}
From the solution of \eqref{eq:IP_system}, we can reconstruct the increment for $\rho$ using the formula
\begin{equation}
\label{eq:Delta_rho}
\Delta \rho = - \left( \URho^{-1} \UNu + \ORho^{-1} \ONu \right) ^{-1}
\left( \res_3 + \URho^{-1} \res_4 - \ORho^{-1} \res_5 - B(u)^\top \Delta u - \Delta \alpha \, e \right) \, .
\end{equation}
The increments for the Lagrange multipliers $\unu$ and $\onu$ are computed based on the stable reduction proposed in \cite{Freund_1997}, with a slight adjustment to account for the upper bound constraints not present in that paper. The multipliers are updated by the following formulas, in the following order
\begin{align}
\label{eq:Delta_onu}
\Delta \onu &=
 \dfrac{1}{\orho-\urho} \left(
  \URho( B(u)^\top \Delta u + \Delta \alpha \, e) - (\UNu - \ONu)\rho - \left( \res_4 + \res_5 - \URho \res_3 \right)
  \right) \, , \\
  \label{eq:Delta_unu}
\Delta \unu &=
 \Delta \onu - B(u)^\top \Delta u - \Delta \alpha \, e - \res_3 \, . 
\end{align}

Once the increments have been obtained, we need to determine an appropriate step length. Our algorithm employs a long step strategy \cite{Wright_1997:IP} in that it restricts the step length mainly to guarantee feasibility of the next iterate. We do not use the same step length for all increments. Rather, $\Delta \rho$ and $\Delta u$ use the same step length, the step length for $\Delta \alpha$ is always equal to 1 and different step lengths are calculated for both $\Delta \unu$ and $\Delta \onu$. For details, see Algorithm \ref{alg:ip}. This strategy proved to be the most effective in numerical experiments.

After each IP iteration, the barrier parameters are updated adaptively. For this, we compute the duality measure for the lower and upper bound constraint
\[
	\dfrac{ \unu^\top (\rho-\urho) }{m} \quad \text{and} \quad \dfrac{ \onu^\top (\orho-\rho) }{m} \,,
\]
respectively. We then scale these measures by a constant $0<\sigma_r<1$ and $0<\sigma_s<1$ to update $r$ and $s$. At this point, one unconventional feature of our algorithm should be highlighted. The new values for $r$ and $s$ are not used to construct the right hand side term for the next iteration, but rather for the iteration after that. We found that this ``iteration shift'', peculiar though it might seem, makes the algorithm significantly more efficient. Indeed, without this shift this version of the code is hardly viable and one requires several Newton iterations per IP iteration instead of just one.

Finally, we require a stopping criterion for the algorithm. Just like in Algorithm~\ref{alg:pbm}, we use the duality gap $\delta(u,\alpha)$ as a measure of optimality, scaled by the current objective function---the \emph{primal} objective function $\frac{1}{2} f^\top u$, in this case. On top of this, we want to ensure that our solution is feasible to within a certain accuracy. Our feasibility measure is the following sum of weighted residuum norms
\begin{equation}
\label{eq:ip_res}
\resIP =
\dfrac{ \| \res_1 \|_2 }{ \| f \|_2 } +
\dfrac{ | \res_2 | }{ V } +
\dfrac{ \| \res_3 \|_2 }{ \| \unu \|_2 + \| \onu \|_2 } +
\dfrac{ | e^\top \res_4 | }{ m } +
\dfrac{ | e^\top \res_5 | }{ m } \, .
\end{equation}
Furthermore, the duality gap should be (nearly) positive, as a negative duality gap points to infeasibility.

Algorithm \ref{alg:ip} sums up our IP method. The parameter values that we used in our experiments are $\epsIP = 10^{-5}$, $\sigma_r=\sigma_s=0.2$. For the initial values, we chose $u=0$, $\alpha=1$, $\rho_i=V/m$ for all $i=1,\dots,m$ and $\unu=\onu=e$. The barrier parameters start at $r=s=10^{-2}$.

\begin{algorithm}
\caption{Primal-dual IP}
\label{alg:ip}
Let $\epsIP>0$ and $0<\sigma_r,\sigma_s<1$ be given. Choose initial vectors $(u,\rho)$ and $(\alpha,\unu,\onu)$. Set barrier parameter update values as $r^+ = \sigma_r \cdot \unu^\top (\rho-\urho)/m$ and $s^+ = \sigma_s \cdot \onu^\top (\orho-\rho)/m$.
\begin{algorithmic}[1]
\Repeat
    \State Solve system \eqref{eq:IP_system} to obtain $(\Delta u, \Delta \alpha)$
    \State Reconstruct $(\Delta \rho, \Delta \onu, \Delta \unu)$ using \cref{eq:Delta_rho,eq:Delta_onu,eq:Delta_unu}
    \State Update barrier parameters:\ 
    $ r = r^+\,,\  s=s^+ $
    \State Compute the following step lengths
    \begin{gather*}
        \kappa_{u}=\kappa_{\rho}= \min \left\{
        0.9 \cdot \min_{\Delta \rho_i>0}{\dfrac{\orho_i-\rho_i}{\Delta\rho_i}} \,,\, 
        0.9 \cdot \min_{\Delta \rho_i<0}{\dfrac{\urho_i-\rho_i}{\Delta\rho_i}} \,, \;
        1 \right\} \\[0.5em]
        \kappa_{\unu} = 0.9 \cdot \min_{\Delta\unu<0} \dfrac{-\unu}{\Delta\unu}, \quad
        \kappa_{\onu} = 0.9 \cdot \min_{\Delta\onu<0} \dfrac{-\onu}{\Delta\onu}, \quad
        \kappa_{\alpha} = 1
    \end{gather*}
    \State Update all variables
    	\begin{gather*}
	u = u + \kappa_u \Delta u \,, \quad
	\alpha = \alpha + \kappa_\alpha \Delta \alpha \,, \quad
	\rho = \rho + \kappa_{\rho} \Delta \rho \,, \\[0.5em]
	\unu = \unu + \kappa_{\unu} \Delta \unu \,, \quad
	\onu = \onu + \kappa_{\onu} \Delta \onu
	\end{gather*}
    \State Compute the duality gap $\delta(u,\alpha)$ by \eqref{eq:gap}, the objective function $\frac{1}{2}f^\top u$ and the feasibility measure $\resIP$ by \eqref{eq:ip_res}
    \If{ $ \epsIP > \delta(u,\alpha) / (\frac{1}{2}f^\top u) > - 0.1\cdot\epsIP$ and $\resIP < 10\cdot\epsIP$}
        \State STOP
    \EndIf
    \State Determine barrier parameters for shifted barrier parameter update
    \[
    r^+ = \sigma_r \cdot \dfrac{ \unu^\top (\rho-\urho) }{m} \,, \quad s^+ = \sigma_s \cdot \dfrac{ \onu^\top (\orho-\rho) }{m}
    \]
\Until convergence
\end{algorithmic}
\end{algorithm}

\section{Optimality Condition (OC) method}\label{sec:OC}
To get a broader picture, we will compare the PBM and IP algorithms with the
established and commonly used Optimality Condition (OC) method. We will
therefore briefly introduce the OC algorithm for VTS. 
For more details, see \cite[p.308]{bendsoe-sigmund} and the references therein. 

We adapt the algorithm implemented in the popular code {\tt top88.m} \cite{top88}; see Algorithm~\ref{alg:doc}. We call it damped OC (DOC) method, due to the exponent $q\leq 1$ that shortens the ``full" OC step. We use the standard value $q=0.5$.
\begin{algorithm}
	\caption{DOC}	
	\label{alg:doc}
	Let $\rho\in\RR^m$ be given such that $\sum_{i=1}^m \rho_i = V$, $ \urho_i \leq \rho_i \leq \orho_i$, $i=1,\ldots,m$. Set $\tau_{\scriptscriptstyle \alpha}=0.1\,\epsDOC$ and $q\leq 1$.
	\begin{algorithmic}[1]
		\Repeat
		\State{$u=(K(\rho))^{-1} f$}
		\State{$\oal=10000$, $\ual=0$}
		\While{$\dfrac{\oal-\ual}{\oal+\ual}>\tau_{\scriptscriptstyle \alpha}$}
		\State{$\alpha = (\oal+\ual)/2$}
		\State{	 $\rho_i^+ = \min\left\{\max\left\{\rho_i \displaystyle\frac{(u^T K_i u)^q}{\alpha}, \urho_i\right\},\orho_i\right\}\,, \quad i=1,\ldots,m$}
		\State{If $\sum_{i=1}^m \rho_i^+>V$ then set $\ual=\alpha$; else if $\sum_{i=1}^m \rho_i^+\leq V$ then set $\oal=\alpha$}
		\EndWhile
        \If{$\|\rho^+-\rho\|_{\inf} \leq \epsDOC$}
        \State STOP
        \EndIf
		\State{$\rho = \rho^+$}
		\Until convergence
	\end{algorithmic}
\end{algorithm}


Following \cite{top88}, we use the stopping criterion
$$
  \|\rho^+-\rho\|_{\inf} \leq \epsDOC\,,
$$
where $\rho$ and $\rho^+$ are the two most recent iterates. While {\tt top88.m}  uses $\epsDOC = 10^{-2}$, we found that this value is too generous in many 3D examples, resulting in an image that is significantly different from an image obtained with $\epsDOC \leq 10^{-3}$; see Figure~\ref{fig:OC} in Section~\ref{sec:Num}, where we address the choice of $\epsDOC$ in a bit more detail.

Another parameter we changed, as compared to  \cite{top88}, was the value of the stopping criterion for the bisection method $\tau_{\scriptscriptstyle \alpha}$. In Section \ref{sec:Num}, we use $\tau_{\scriptscriptstyle \alpha}=0.1\,\epsDOC$ which leads to a more stable behaviour of the DOC and only marginal increase of total CPU time.

The reader may ask about the relation of the DOC stopping criterion (using difference of variables in two subsequent iterations) with the more rigorous criterion based on the duality gap, used in the PBM and IP algorithms. Our experiments revealed a somewhat surprising phenomenon: in most of the problems we solved, the behaviour of the two stopping measures was almost identical. This experience justifies the use of the DOC stopping criterion and, in particular, the relative fairness of our comparisons of DOC with PBM and IP. 
\section{Multigrid preconditioned Krylov subspace methods}\label{sec:MG}
In the previous sections, we have introduced three algorithms for the solution of the VTS problem, all of which have one thing in common: In every iteration, we have to solve a system of linear equations
\begin{equation}\label{eq:lineq}
Az=b \,,
\end{equation}
where $b\in\RR^n$ and $A$ is a $n\times n$ symmetric positive definite matrix.
In the OC method, $A$ is the stiffness matrix $K(\rho)$  of the linear elasticity problem. In algorithms PBM and IP, $A$ corresponds to the Schur complements $S$ from equations \eqref{eq:S} and \eqref{eq:S_IP}, respectively. These latter two matrices have the same sparsity structure. In particular, the principal $(n-1)\times (n-1)$ submatrix of $S$ has the same sparsity structure as the stiffness matrix $K(\rho)$; the last row and column of $S$ are full.

In this section, we will recall an iterative method that is known to be very efficient for linear elasticity problems on well structured finite element meshes.
Throughout the section, we will use the notation of \eqref{eq:lineq}.

\subsection{Multigrid preconditioned MINRES}
We use standard V-cycle correction scheme multigrid method with coarse level problems
$$
A_k z^{(k)}=b^{(k)},\quad k=1,\ldots,\ell-1\,,
$$
where
$$
A_{k-1} = I_k^{k-1} (A_k) I_{k-1}^k,\quad b^{(k-1)} = I_k^{k-1} (b^{(k)}),\quad k=2,\ldots,\ell\,.
$$
Here we assume that there exist $\ell-1$ linear operators
$I_k^{k-1}:\RR^{n_k}\to\RR^{n_{k-1}}$, $k=2,\ldots,\ell$, with
$n:=n_\ell>n_{\ell-1}>\cdots>n_2>n_1$ and let $I_{k-1}^k:=(I_k^{k-1})^T$.
As a smoother, we use the Gauss-Seidel iterative method.
See, e.g., \cite{hackbusch} for details.

Although the multigrid method is very efficient, an even more
efficient tool for solving  (\ref{eq:lineq}) may be a preconditioned
Krylov type method, whereas the preconditioner consists of one V-cycle of the multigrid method\footnote{We found more than one V-cycle to not be as efficient in terms of overall CPU time.}. After experimenting with several Krylov methods, we found that the MINRES algorithm \cite{paige1975solution} is the most robust for our problems in which the system matrix may converge to a positive semidefinite matrix. We use the standard implementation of MINRES from \cite{barrett1994templates}.

\subsection{Multigrid MINRES for PBM, IP and OC}

In all examples in Section~\ref{sec:Num}, we use hexahedral finite elements with trilinear basis functions for the displacement variable $u$ and constant basis functions for the
variable $\rho$, as is the standard in topology optimization. We start with a very coarse mesh and use regular refinement of each element into 8 new elements. The
prolongation operators $I_{k-1}^k$ for the variable $u$ are based on a
standard 27-point interpolation scheme.
For more details, see, e.g., \cite{hackbusch}.
When solving the linear systems \eqref{eq:system} and \eqref{eq:IP_system} in PBM and IP,
we also need to prolong and restrict the single additional variable $\lambda$;
here we simply use the identity.

When we use the regular finite element refinement mentioned above, the stiffness matrix $K(\rho)$ will be sparse and, if a reasonably good numbering of the nodes is used, banded. The number of non-zero elements in a row of $K(\rho)$ does not exceed 81. A typical non-zero structure of $K$ is
shown in Figure~\ref{subfig:1right}, if we ignore the additional last
column and row in that figure.

As usual, the MINRES method is stopped whenever
\begin{equation}\label{eq:cgstop}
\| r \|\leq \epsMINRES\|b\| \,,
\end{equation}
where $ r $ is the residuum and $b$ the right-hand side of the linear system,
respectively. The choice of the stopping parameter $\epsMINRES$ varies between the different algorithms.

\paragraph{Multigrid MINRES for OC}\label{sec:CGOC}
The only degree of freedom in the algorithm is the stopping
criterion. The required accuracy of these solutions (such that the
overall convergence is maintained) is well documented and theoretically
supported in the case of the IP method; it is, however, an unknown in the case of the DOC
method; see \cite{amir} for detailed discussion. Clearly, if the linear systems
in the DOC method are solved too inaccurately, the whole method may diverge or
just oscillate around a point which is not the solution.

In all our numerical experiments, we used $\epsMINRES=10^{-4}$.
In \cite{MK_Mohammed_2016}, it was observed that, with this stopping criterion, the number of
DOC iterations was almost always the same, whether we used an iterative or a
direct solver for the linear systems. Our experiments with 3D problems confirmed this observation.

\paragraph{Multigrid MINRES for PBM}
The initial stopping parameter $\epsMINRES$ scales with the size of the problem, as it can otherwise be too strict for large problems or too imprecise for small problems. We initialize and update it in the following way:
\begin{itemize}
	\item We start with $\epsMINRES=10^{-4}\sqrt{n}$.
	\item Let $\resPBM$ be the sum of the residua computed in the current step of the PBM Newton Algorithm~\ref{alg:pbm_newton} and let $\resPBM^+$ be this sum in the following step. If $\resPBM^+>0.9\,\resPBM$, we update $\epsMINRES := \max\{0.1\,\epsMINRES,10^{-9}\}$. In other words, we increase the accuracy of the stopping parameter whenever we do not achieve a satisfactory improvement in feasibility and optimality with the current $\epsMINRES$.
\end{itemize}

In our numerical tests, the update had to be done only in a few cases and the
smallest value of $\epsMINRES$ needed was  $\epsMINRES=10^{-3}$.

\paragraph{Multigrid MINRES for IP}
In the IP method, we use an adaptive updating scheme for the stopping parameter, based on the complementarity of the current solution:
\begin{itemize}
	\item We start with $\epsMINRES = 10^{-2}$
	\item We compute
		$$ d = \max \left\{
			\max_{i=1,\dots,m} | (\rho_i - \urho_i) \unu_i | \,, \max_{i=1,\dots,m} | (\orho_i - \rho_i) \onu_i |
			\right\}
			$$
		and set $\epsMINRES := \max\{ 100\,d, 10^{-9} \} $ if this new value is lower than the current $\epsMINRES$. The low minimum value of $10^{-9}$ for $\epsMINRES$ has proven to be necessary for convergence in our experiments. 
\end{itemize}

\section{Numerical experiments}\label{sec:Num}
We now present and compare numerical results for the PBM, IP and DOC methods. In Section \ref{subsec:m-scale}, we focus on a rigorous comparison of the performance of the three algorithms, both in terms of CPU time and required calls to the iterative solver. For this, we look at problems where the number of finite elements is in the order of $10^4$ to $10^5$. As we will see, the PBM method outperforms both the IP and the DOC methods. When we consider problems with over a million finite elements in Section \ref{subsec:l-scale}, we only present results for IP and PBM, since DOC with our required accuracy is no longer practicable.

In the formulation of the VTS problem \eqref{eq:to_intro}, we chose the lower bounds $\urho$ to be positive. As far as the underlying physical model is concerned, however, $\urho=0$ would make the most sense, with $\rho_i=0$ corresponding to an element without material. A lower bound larger than zero might distort the, as it were, physically more accurate results. Yet the strict positivity is required for the positive definiteness of $K(\rho)$ and to bound the condition number of the system matrices arising in the different methods. In our experiments, this turned out to be critical for the OC and IP, but not for the PBM method. Therefore, we generally set $\urho=0$ for PBM only.

The code was implemented in Matlab, outsourcing certain subroutines to C via MEX files. No parallelization was performed in any of our functions. While the Matlab inbuilt routines are in general parallelizable, on the BlueBEAR HPC system used to produce the large-scale results in Section~\ref{subsec:l-scale}, it was limited to a single core.

The design domain for each of our example problems is set up in a way that is based on a multigrid-structure. It is a cuboid defined by $m_x\times m_y\times m_z$ cubes of equal size corresponding to the coarse level finite element mesh. 
We refine the coarse mesh regularly $\ell-1$ times, giving us $\ell$ mesh levels in total; each cube element is refined into 8 new elements of equal dimensions. Hence level-2 refinement of a $4\times 2 \times 2$ coarse mesh with 16 elements results in a $8\times 4 \times 4$ mesh with 128 element and level-$\ell$ refinement of the same coarse mesh results in a mesh with $16\cdot 8^{\ell-1}$ elements.

\begin{figure}[h]
	\begin{center}
		\begin{subfigure}[b]{0.48\textwidth}
			\includegraphics[width=\hsize]{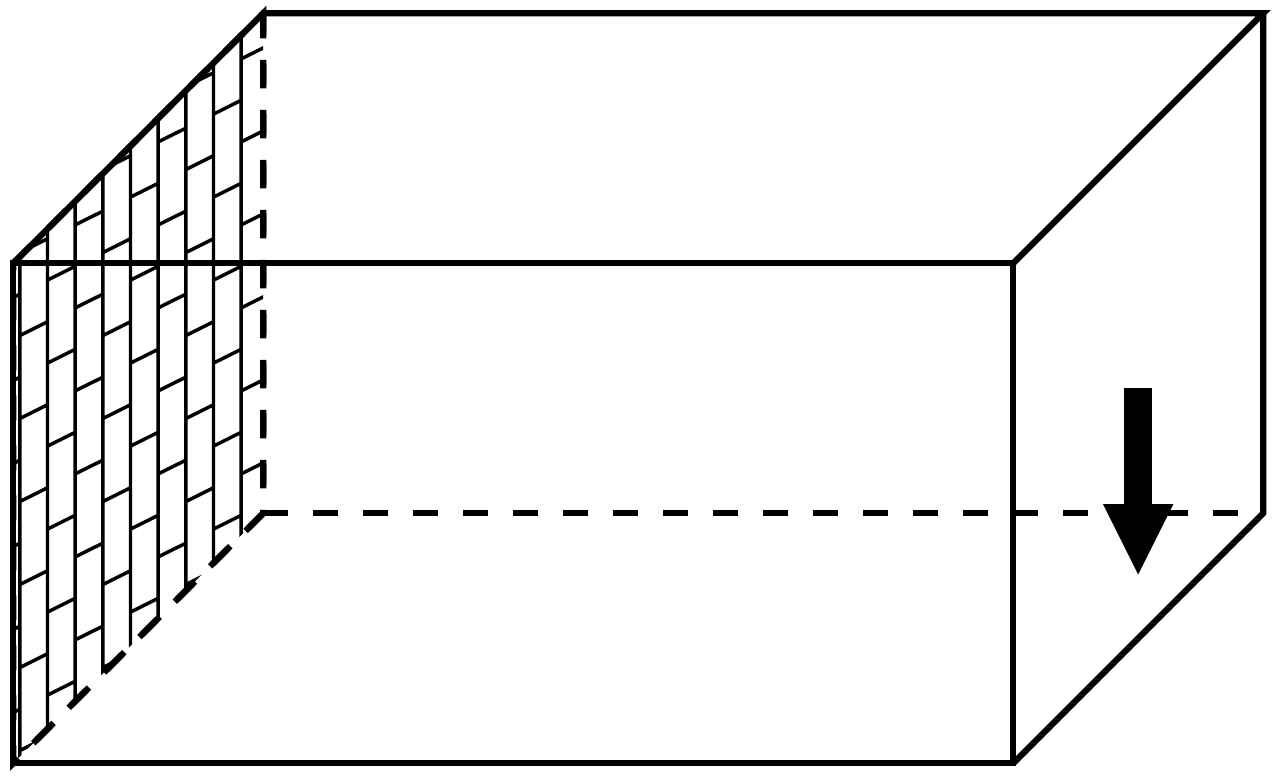}
			\caption{cantilever}
			\label{subfig:mesh2}
		\end{subfigure} \quad
		\begin{subfigure}[b]{0.48\textwidth}
			\includegraphics[width=\hsize]{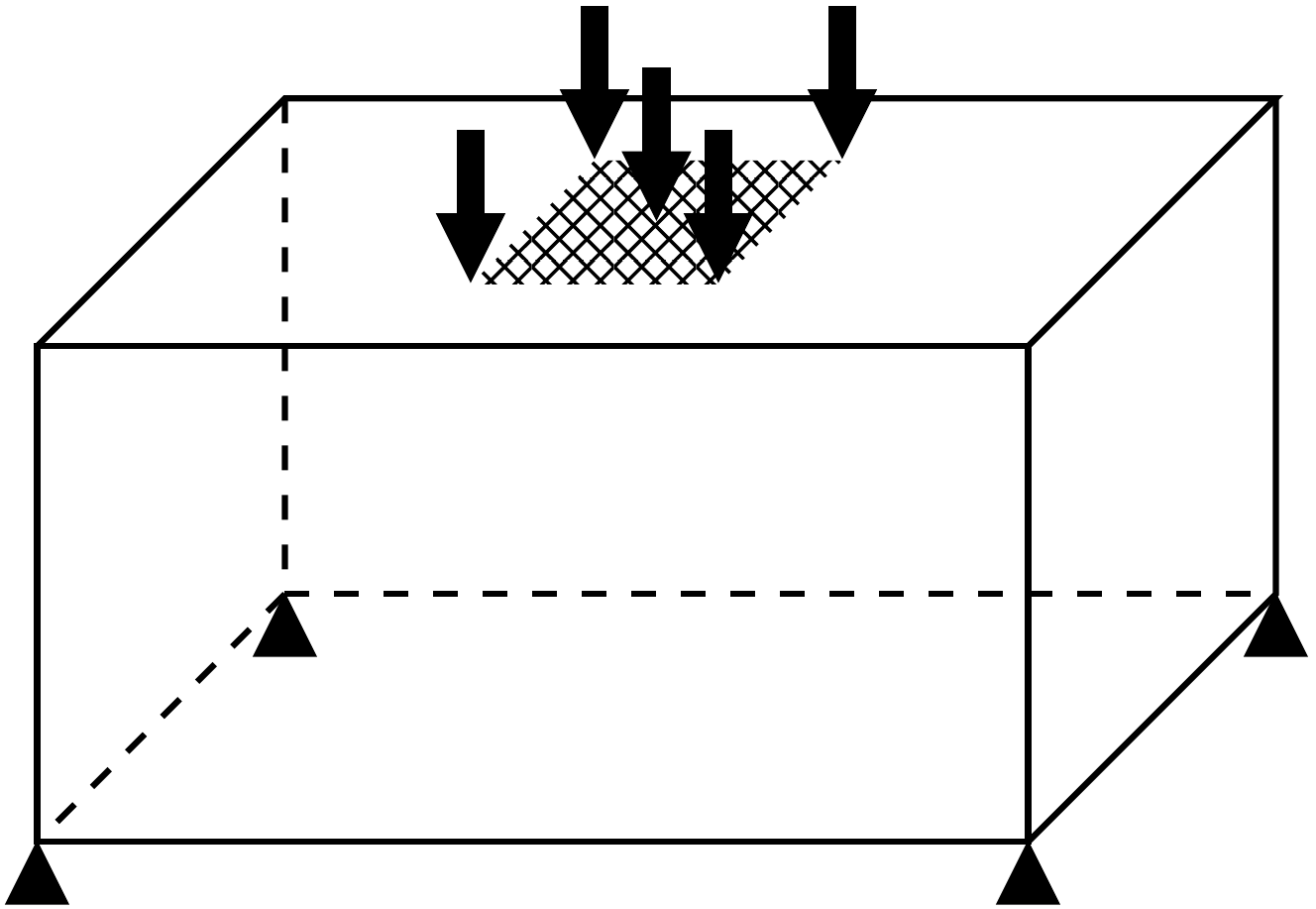}
			\caption{bridge}
			\label{subfig:mesh3}
		\end{subfigure}
	\end{center}	
	\caption{\label{fig:mesh23}Boundary conditions and loads for cantilever and bridge problems.}
\end{figure}

We consider two sets of boundary conditions and loading scenarios, referring to the first one as ``cantilever" and to the second one as ``bridge"; see Figure~\ref{fig:mesh23} for the specifications. Cantilever problems have all nodes on the left-hand face fixed in all directions; a load in direction $z$ is applied in the middle of the right-hand face (Figure~\ref{subfig:mesh2}). The bridge problems are subject to a uniform load applied on a rectangle centered on the top face; all four corners of the bottom face are fixed in all directions (Figure~\ref{subfig:mesh3}). 

We adopt the following naming convention for the problems solved in this paper:
\begin{description}
	\item{\bfseries CANT-{\boldmath $m_x$}-{\boldmath $m_y$}-{\boldmath $m_z$}-{\boldmath $\ell$}} for a cantilever with a $m_x\times m_y\times m_z$ coarse mesh and $\ell$ mesh levels;
	\item{\bf BRIDGE-{\boldmath $m_x$}-{\boldmath $m_y$}-{\boldmath $m_z$}-{\boldmath $\ell$}} for a bridge with a $m_x\times m_y\times m_z$ coarse mesh and $\ell$ mesh levels.
\end{description} 

\subsection{Comparison of PBM, IP and OC}\label{subsec:m-scale}
The problems in this section have been solved on a 
2018 MacBook Pro with 2.3GHz dual-core Intel Core i5, Turbo Boost up to 3.6GHz, and 16GB RAM.
This allowed us to properly compare the CPU timing; but it also prevented us from solving large scale problems, due to memory limitations. The results for those problems, run on a HPC computer, are reported in the next section.


\paragraph{Example CANT-16-2-2-5}
In Table~\ref{tab:2} we present results for problem CANT-16-2-2-5 with 262\,144 finite elements. The lower bound for $\rho$ was set to zero for the PBM method and to $\urho=10^{-7}$ for the IP and DOC methods. 

Each table row shows the results for a certain method and stopping parameter. They are given in terms of the total number of linear systems solved\footnote{This is equal to the number of Newton iterations in the case of the PBM method. For the other two algorithms, there is no difference between the number of ``outer'' iterations and the number of Newton iterations.}; the total number of MINRES iterations; the total CPU time needed to solve the problem; the CPU time spent on solving the linear systems; and the final value of the primal objective function, where the accurate digits\footnote{Digits are assumed to be accurate when the different methods all appear to converge to them.} are in bold.

Because the IP method had difficulties getting below our stopping threshold $\epsIP=10^{-5}$, we also ran this method with $\urho=10^{-3}$ for comparison, since this improves the conditioning of the system matrix. The resulting objective value is not comparable with the other values and is thus grayed out.

\begin{figure}[h]
	\begin{center}
		\begin{subfigure}{0.78\textwidth}
			\includegraphics[width=\hsize]{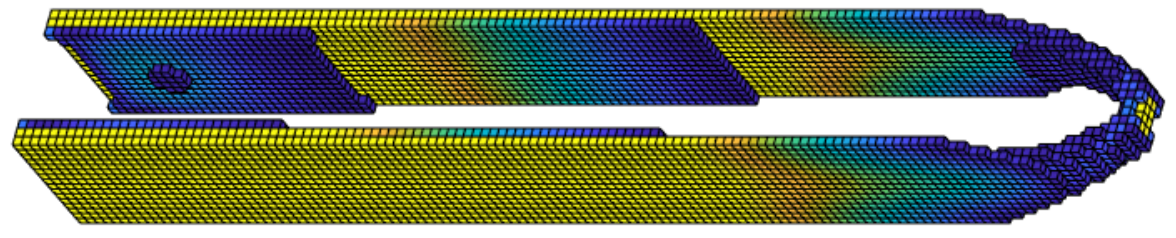}
			\caption{$\epsDOC = 10^{-2}$}
			\label{subfig:OCa}
		\end{subfigure} \quad
			\begin{subfigure}{0.78\textwidth}
		\includegraphics[width=\hsize]{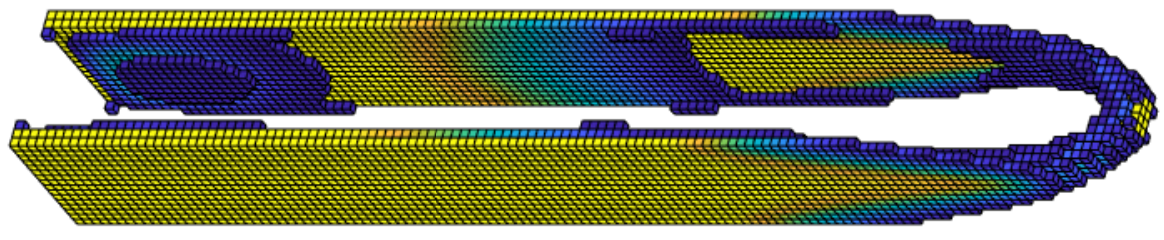}
		\caption{$\epsDOC = 10^{-3}$}
		\label{subfig:OCb}
	\end{subfigure} \quad
		\begin{subfigure}{0.78\textwidth}
			\includegraphics[width=\hsize]{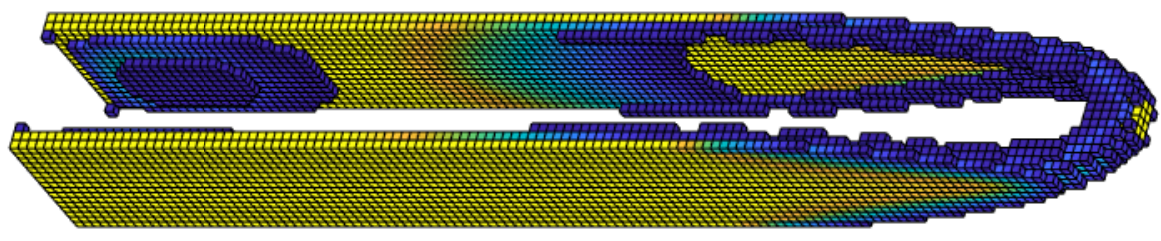}
			\caption{$\epsDOC = 10^{-5}$}
			\label{subfig:OCc}
		\end{subfigure}
	\end{center}	
	\caption{\label{fig:OC}CANT-16-2-2-4, DOC result with $\epsDOC = 10^{-2}$, $\epsDOC = 10^{-3}$, $\epsDOC = 10^{-5}$. Figure (c) is identical with IP and PBM results. Only elements with density values of $\rho_i>0.1$ are shown in order to make the differences visible.}
\end{figure}

We ran the DOC method with three different stopping tolerances. While $\epsDOC=10^{-2}$ would be used to mimic the {\tt top88} code, we can see in Figure~\ref{fig:OC} that the final result delivered with this tolerance is by no means optimal and clearly differs from that obtained with $\epsDOC=10^{-5}$ (for better transparency, Figure~\ref{fig:OC} present results of a smaller problem CANT-16-2-2-4). Decreasing $\epsDOC$ to $10^{-3}$ improves the result but the image is still visibly different from the optimal one. This is despite the five correct significant digits in the objective function, reached by DOC with $\epsDOC=10^{-3}$.
The results produced by PBM and IP were ``visually identical" to that for DOC with $\epsDOC=10^{-5}$ in Figure~\ref{subfig:OCc}. (Of course, this ``visual comparison" is not rigorous but, in the end, the image is the required result of topology optimization in practice; a rigorous comparison is given in Table~\ref{tab:2}.)

We also ran the PBM method with a lower stopping tolerance $\epsPBM=10^{-6}$ to demonstrate that the method can reach higher precision with only relatively few additional iterations.

The numbers in Table~\ref{tab:2} show that PBM clearly outperforms the other two methods, both with respect to the number of MINRES iterations and to the CPU time required by the whole algorithm and the linear solver only. It is even faster than the DOC method with the very relaxed stopping tolerance $\epsDOC=10^{-2}$, at the same time delivering a solution of much higher quality.
\begin{table}[htbp]
	\centering
	\caption{Example CANT-16-2-2-5 by different methods. Problem dimensions: $m = 262\,144,\ n=836\,352$.}
	\renewcommand{\textbf}[1]{\fontseries{b}\selectfont #1\normalfont}
	\begin{tabular}{l r r r r r S[table-number-alignment = center, detect-weight, mode=text, table-format = 3.6]}
		\toprule
		& stop&\multicolumn{2}{c}{iterations} & \multicolumn{2}{c}{CPU time [s]} \\
		method &tol  & Nwt/OC  & MINRES  & total & lin~solv & \multicolumn{1}{r}{obj fun}\\
		\midrule
		PBM   & $10^{-5}$ & 42      &     156  &   916  &  317 & \textbf{66.192}8136\\
		PBM   & $10^{-6}$ & 48      &     259   &  1110  &  420 & \textbf{66.19273}18\\
		IP       & $10^{-5}$ & 36      &   4992   & 6510  & 5670 & \textbf{66.19272}49\\
		IP($\urho=10^{-3}$)  & $10^{-5}$   & 29      &    1243   & 1510   &  1070 &  {\color{gray} 66.1988863}\\
		DOC   & $10^{-2}$ & 38      &    394   & 1160   &  883 &  \textbf{66.2}107223\\
		DOC   & $10^{-3}$ & 226     &  2462    & 6710  & 5020 &  \textbf{66.193}4912\\
		DOC   & $10^{-5}$ & 2759   &  30325 & 82500  & 61900 & \textbf{66.19272}72\\
		\bottomrule
	\end{tabular}%
	\label{tab:2}%
\end{table}%

Below are some further, detailed observations:
\begin{itemize}
	\item The PBM iterations are very robust in terms of MINRES iterations needed to solve the linear systems. Up to the very last PBM iterations, MINRES only requires 1--3 steps to reach the required accuracy. Even in the last PBM iterations, the number of MINRES steps typically does not exceed 15--20. One reason for this is presumably that the updating scheme for the MINRES tolerance $\epsMINRES$ (see Section~\ref{sec:MG}) only rarely needs to update the value. With $\epsMINRES$ thus decreasing only very slowly, the linear systems never have to be solved to a very high accuracy. Still, the PBM solution displays the required optimality and feasibility.
	\item The IP method is much more sensitive to ill-conditioning. While in the first IP iterations MINRES only requires 1--2 steps, this number then quickly increases when nearing the required IP stopping criterion. In the CANT-16-2-2-5 problem with $\urho=10^{-3}$, the number of MINRES steps in the IP Newton iterations grew as follows: 1--1--1--1--1--1--1--1--1--1--1--2--3--3--6--5--7--11--13--23--35--49--55--25--66-466--149--314. 
	\item The number of MINRES steps in every DOC iteration is almost constant. In the CANT-16-2-2-4 problem, this number was between 8 and 11 in the first 49 DOC iterations and 12 for all remaining DOC iterations, even with the stopping tolerance $\epsDOC=10^{-5}$. The total number of DOC iterations, however, grows dramatically when higher precision in the stopping criterion is required.
    \item Because of the way that $\rho$ is computed in the different algorithms, the volume constraint is not satisfied to the same degree of accuracy in each case. The OC method yields the most accurate $\rho$ with respect to the volume constraint, while the PBM solution generally gives $\sum_i \rho_i > V$. The PBM solution deviation from $V$ was never more than one permille in our experiments.
\end{itemize}

\paragraph{Example BRIDGE-4-2-2-6}
We now present some results of the BRIDGE problem. Table~\ref{tab:5} shows the iteration numbers and CPU times for BRIDGE-4-2-2-6 with 524\,288 finite elements. Compared to CANT-16-2-2-x, the stiffness matrix in these problems (and thus the Schur complement for each method) has a higher condition number, due to the different shape of the computational domain. 
\begin{table}[htbp]
	\centering
	\caption{Example BRIDGE-4-2-2-6 solved by different methods. Problem dimensions: $m = 524\,288,\ n=1\,635\,063$.}
	\renewcommand{\textbf}[1]{\fontseries{b}\selectfont #1\normalfont}
	\begin{tabular}{l r r r r r S[table-number-alignment = center, detect-weight, mode=text, table-format = 3.6]}
		\toprule
		& stop&\multicolumn{2}{c}{iterations} & \multicolumn{2}{c}{CPU time [s]} \\
		method &tol  & Nwt/OC  & MINRES  & total & lin~solv & \multicolumn{1}{r}{obj fun}\\
		\midrule
		PBM   & $10^{-5}$ & 57      &     330  &  2710  &  1020 & \textbf{42.000}2293\\
		PBM   & $10^{-6}$ & 62      &     423  &  3000  &  1190 & \textbf{42.00015}61\\
		IP    & $10^{-5}$ & 49      &    2919  &  6040  &  4320 & \textbf{42.000}2076\\
		IP($\urho=10^{-3}$)  & $10^{-5}$   & 51      &    2965   & 6210   &  4440 &  {\color{gray} 42.0027639}\\
		DOC   & $10^{-2}$ & 99      &    1454  &  6800  &  5330 & \textbf{42.00}14281\\
		DOC   & $10^{-3}$ & 278     &    4139  & 19200  & 15100 & \textbf{42.000}2175\\
		DOC   & $10^{-5}$ & 659     &    9854  & 45900  & 36100 & \textbf{42.00015}23\\
		\bottomrule
	\end{tabular}%
	\label{tab:5}%
\end{table}%


\subsection{Large scale problems}\label{subsec:l-scale}
In this section, we do not include the CPU times needed to solve the example problems. This is because they were solved on the Linux HPC BlueBEAR with 2000 cores of different types, with up to 498 GB RAM per core. We did not have any control over which cores were used for which job, so that the time statistics could not be used for reliable performance comparison. Furthermore, recall that Matlab only ran on a single core on BlueBEAR, so that the total computation time for any example would most likely not be competitive compared with any parallelized code.

We present results for the PBM and IP algorithms only. As we have seen in Section~\ref{subsec:m-scale}, they are both several times faster than the DOC method for the same degree of accuracy. This does not improve with larger problem sizes, which means that the OC method might take several days to solve a problem which is solved in just a few hours by the PBM method. To solve the BRIDGE-$4$-$2$-$2$-$5$ and BRIDGE-$4$-$2$-$2$-$6$ problems, for example, the OC method requires roughly 17 times as much CPU time as the PBM method. This factor is even larger for CANT-$16$-$2$-$2$-$4$ and CANT-$16$-$2$-$2$-$5$. Comparisons for CANT-$4$-$2$-$2$-$5$ and CANT-$4$-$2$-$2$-$5$, which are not included here, gave a factor of over 20.

As before, we set $\urho=0$ for the PBM method. For IP, we chose $\urho=10^{-3}$, as ill-conditioning becomes critical in the large-scale problems covered in this section. Even with this lower bound, IP was not able to solve all of the examples we considered. When it failed, no convergence was apparent once the duality gap had gotten below a certain threshold, which was typically still two or three orders of magnitude too large for the stopping criterion.

The same two problems are considered as in the previous section, namely CANT-$m_x$-$m_y$-$m_z$-$\ell$ and BRIDGE-$m_x$-$m_y$-$m_z$-$\ell$. In this section, we fix the width and height of the design domain to $m_y=m_z=2$ and vary the length $m_x=2,4,6,8$. We ran the code with $\ell=5,6,7$ mesh levels. Tables~\ref{tab:CANT_large} and \ref{tab:BRIDGE_large} show the results for \mbox{CANT-$m_x$-$m_y$-$m_z$-$\ell$} and BRIDGE-$m_x$-$m_y$-$m_z$-$\ell$ in terms of iteration numbers.

\begin{figure}    
    \caption{Optimal density $\rho$ for CANT-8-2-2-7. The elements with the lowest density values are hidden such that the visible element densities add up to $0.8\cdot V$.}
    \centering
    \includegraphics[trim=0 60pt 20pt 40pt, clip, width=0.8\textwidth]{%
        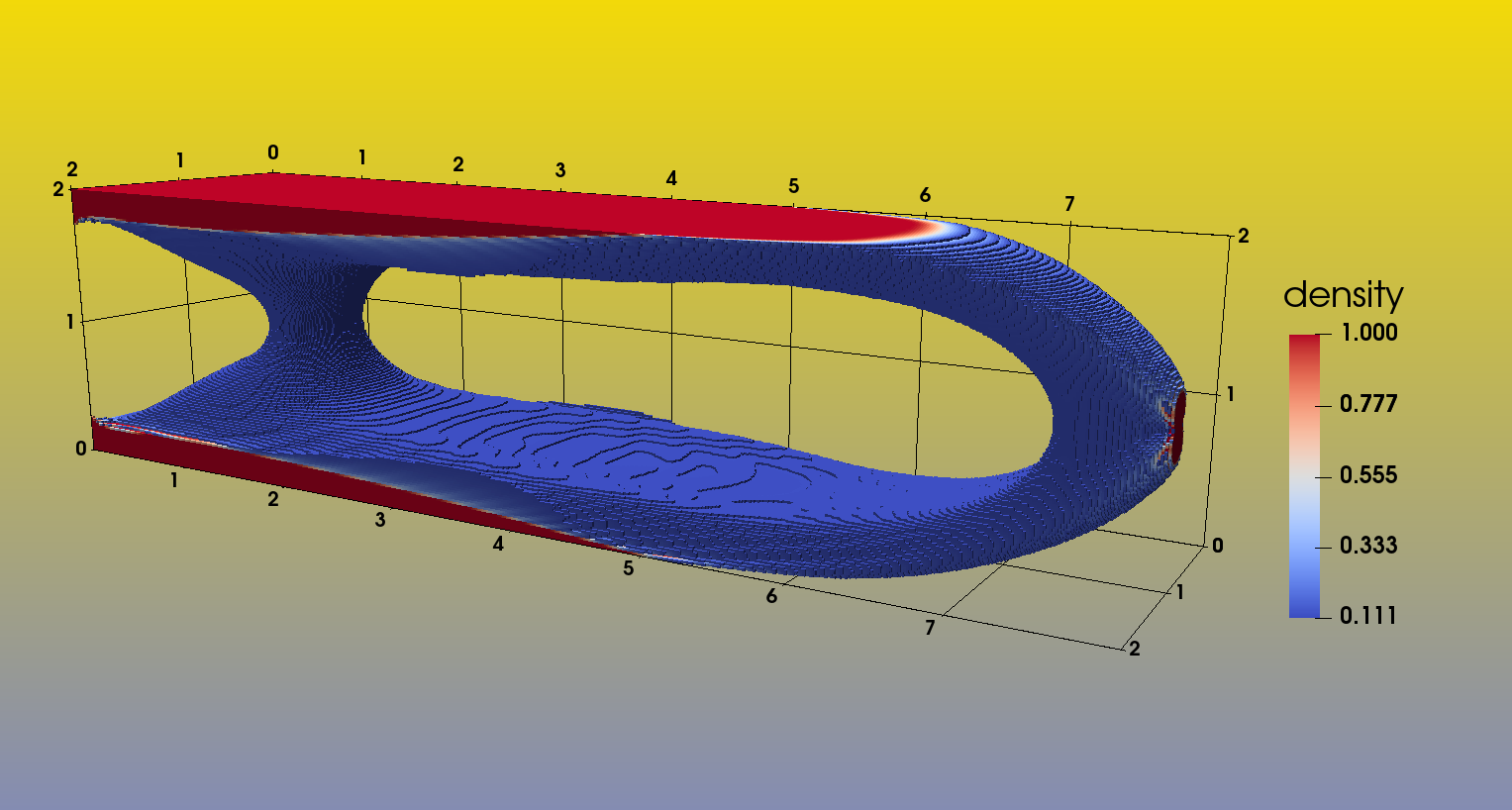}
    \label{fig:CANT8227}
\end{figure}

\begin{figure}  
    \caption{Optimal density $\rho$ for BRIDGE-8-2-2-7. The elements with the lowest density values are hidden such that the visible element densities add up to $0.8\cdot V$.}
    \centering
    \includegraphics[trim=0 60pt 20pt 40pt, clip, width=0.8\textwidth]{%
        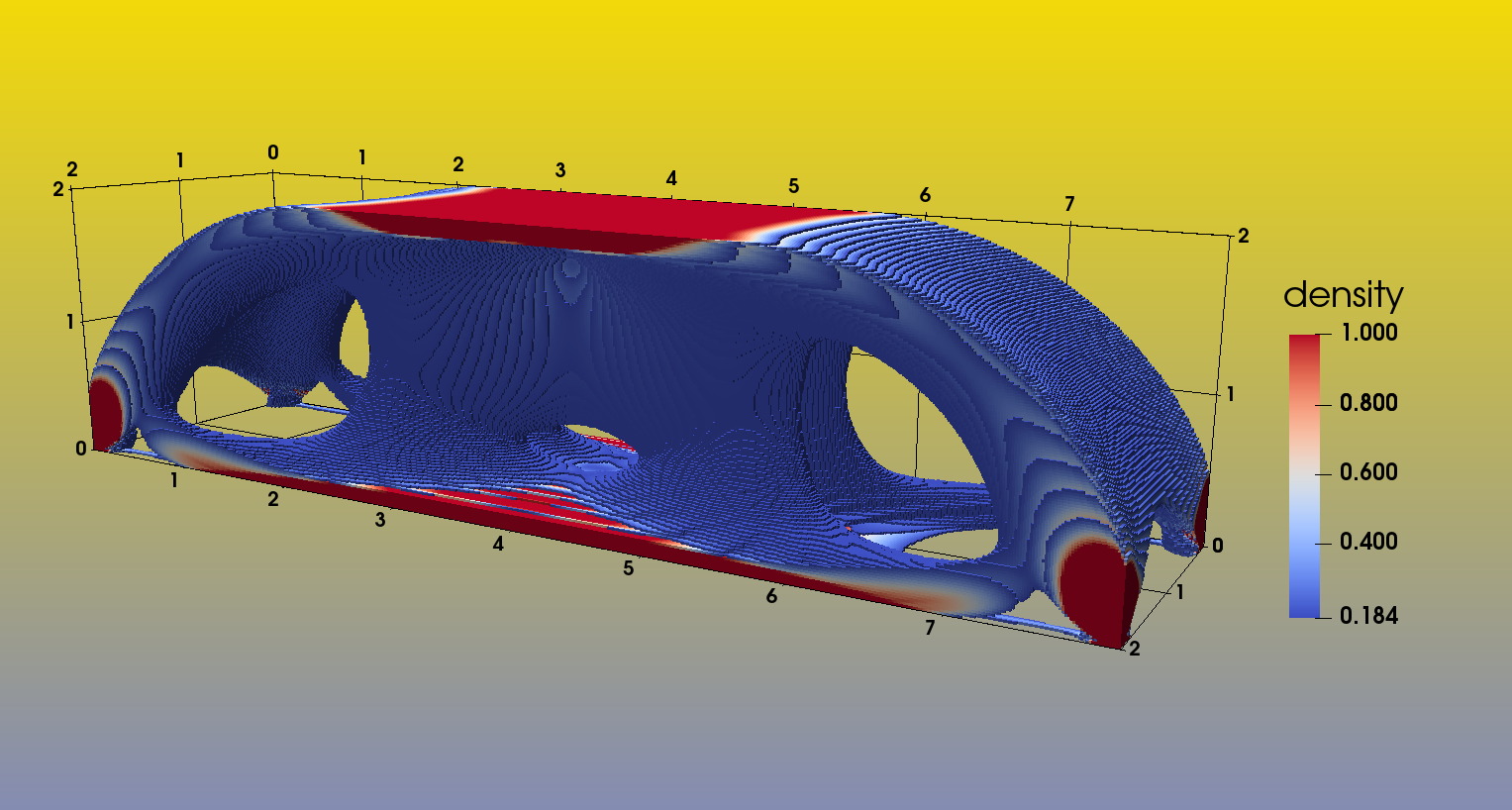} 
    \label{fig:BRIDGE8227}
\end{figure}

The optimal designs produced by the PBM method can be seen in Figures~\ref{fig:CANT8227} and \ref{fig:BRIDGE8227}. The VTS solution typically has a large ``gray area'', i.e., $\rho_i$ is well within the interval $[\,\urho,\orho\,]$ for the majority of elements. This makes it less straightforward to interpret the solution as a discrete design than it is in the case of the SIMP formulation \cite{bendsoe-sigmund}. We must determine a cut-off value $\rho^{*}$ such that all elements with $\rho_i<\rho^{*}$ are ignored. As the design domain is elongated, the density distribution further does not change in a linear fashion. Rather, the gray area is spread disproportionately more thin while most solid elements are clustered along the boundary. Therefore, instead of choosing a constant cut-off value, we found that the most consistent way to plot the results was to consider only the densest elements which add up to a fixed proportion $cV$ of allowed volume, where we chose $c=0.8$.

\newcommand*{\nc}{-- & --}
\begin{table}
    \caption{Example CANT-$m_x$-$m_y$-$m_z$-$\ell$ solved by IP and PBM. Overall IP/PBM iterations, Newton iterations and MINRES iterations. Non-default parameters: (1)  $\gamma=\beta=0.5$ and initial $\epsMINRES=10^{-5}\sqrt{n}$.}
    \label{tab:CANT_large}
    \centering
    \begin{tabular}{l
            @{\hspace{0.8em}}
            S[table-number-alignment = right, table-figures-integer=8,
            table-figures-decimal = 0]
            S[table-number-alignment = right, table-figures-integer=8,
            table-figures-decimal = 0]
            *{5}r }
        \toprule
        \multicolumn{3}{c}{Problem dimensions} & \multicolumn{2}{c}{IP} & \multicolumn{3}{c}{PBM} \\
        \cmidrule(lr){1-3} \cmidrule(lr){4-5} \cmidrule(lr){6-8}
        $m_x$-$m_y$-$m_z$-$\ell$ & 
        \multicolumn{1}{r}{$m$} & \multicolumn{1}{r}{$n$}
         & IP/Nwt & MR & PBM & Nwt & MR \\
        \midrule
        2-2-2-5 & 32768 & 104544 & 31 & 368 & 16 & 50 & 175 \\
        4-2-2-5 & 65536 & 209088 & 28 & 570 & 15 & 57 & 153 \\
        6-2-2-5 & 98304 & 313632 & 26 & 467 & 14 & 45 & 84 \\
        8-2-2-5 & 131072 & 418176 & 27 & 489 & 14 & 45 & 115 \\[0.5em]
        2-2-2-6 & 262144 & 811200 & 46 & 1195 & 18 & 60 & 141 \\
        4-2-2-6 & 524288 & 1622400 & 42 & 2465 & 17 & 59 & 118 \\
        6-2-2-6 & 786432 & 2433600 & 39 & 1015 & 17 & 66 & 157 \\
        8-2-2-6 & 1048576 & 3244800 & 39 & 1079 & 16 & 57 & 88 \\[0.5em]
        2-2-2-7 & 2097152 & 6390144 & 71 & 2383 & 22 & 66 & 70 \\
        4-2-2-7 & 4194304 & 12780288 & 54 & 3543 & 20 & 57 & 67 \\
        6-2-2-7 & 6291456 & 19170432 & 57 & 2667 & 19 & 60 & 68 \\
        8-2-2-7 & 8388608 & 25560576 & 58 & 2335 & 29\rlap{\hspace{1pt}${}^1$} & 64\rlap{\hspace{1pt}${}^1$} & 100\rlap{\hspace{1pt}${}^1$} \\
        \bottomrule
    \end{tabular}
\end{table}

\begin{table}
    \caption{Example BRIDGE-$m_x$-$m_y$-$m_z$-$\ell$ solved by IP and PBM. Overall IP/PBM iterations, Newton iterations and MINRES iterations. When $\resPBM$ in the final PBM iteration did not go below $\epsNWT$, the value at the accepted solution is given. Non-default parameters: (1) $\gamma=\beta=0.5$; (2) initial $\epsMINRES=10^{-5}\sqrt{n}$; (3) initial $\epsMINRES=10^{-5}\sqrt{n}$ and $\epsNWT=0.1$. }
    \label{tab:BRIDGE_large}
    \centering
    \begin{tabular}
            {l @{\hspace{0.8em}}
            S[table-number-alignment = right, table-figures-integer=8,
                    table-figures-decimal = 0]
            S[table-number-alignment = right, table-figures-integer=8,
            table-figures-decimal = 0]
            *{5}r @{\hspace{1.7em}} S[table-format = 1.2e2] 
        }
        \toprule
        \multicolumn{3}{c}{Problem dimensions} & \multicolumn{2}{c}{IP} & \multicolumn{4}{c}{PBM} \\
        \cmidrule(lr){1-3} \cmidrule(lr){4-5} \cmidrule(lr){6-9}
         & \multicolumn{1}{r}{$m$} & \multicolumn{1}{r}{$n$} & IP & MR & PBM & Nwt & MR &
        \multicolumn{1}{c}{$\resPBM$} \\
        \midrule
        2-2-2-5 & 32768 & 107799 & 24 & 387 & 15 & 49 & 220 \\
        4-2-2-5 & 65536 & 212343 & 25 & 590 & 14 & 45 & 155 \\
        6-2-2-5 & 98304 & 316887 & 26 & 778 & 15 & 55 & 309 \\
        8-2-2-5 & 131072 & 421431 & 25 & 1050 & 13 & 47 & 184 \\[0.5em]
        2-2-2-6 & 262144 & 823863 & 41 & 2029 & 15 & 56 & 263 \\
        4-2-2-6 & 524288 & 1635063 & 50 & 2965 & 15 & 57 & 317 \\
        6-2-2-6 & 786432 & 2446263 & \nc & 15 & 61 & 466 \\
        8-2-2-6 & 1048576 & 3257463 & \nc & 15 & 62 & 592 \\[0.5em]
        2-2-2-7 & 2097152 & 6440055 & 91 & 4744 & 26\rlap{\hspace{1pt}${}^1$} & 109\rlap{\hspace{1pt}${}^1$} & 1134\rlap{\hspace{1pt}${}^1$} & 1.14e-4 \\
        4-2-2-7 & 4194304 & 12830199 & \nc & 26\rlap{\hspace{1pt}${}^1$} & 99\rlap{\hspace{1pt}${}^1$} & 718\rlap{\hspace{1pt}${}^1$} & 2.68e-4 \\
        6-2-2-7 & 6291456 & 19220343 & \nc & 25\rlap{\hspace{1pt}${}^{1,3}$} & 98\rlap{\hspace{1pt}${}^{1,3}$} & 743\rlap{\hspace{1pt}${}^{1,3}$} & 2.22e-4 \\
        8-2-2-7 & 8388608 & 25610487 & \nc & 25\rlap{\hspace{1pt}${}^{1,2}$} & 97\rlap{\hspace{1pt}${}^{1,2}$} & 707\rlap{\hspace{1pt}${}^{1,2}$} & 1.17e-3  \\
        \bottomrule
    \end{tabular}
\end{table}

To solve some of the examples by the PBM method, we had to deviate from the choice of parameters specified earlier. For some examples with $\ell=7$ refinement levels, we set $\gamma=\beta=0.5$, rather than $\gamma=\beta=0.3$. Otherwise, the penalty parameters are scaled down too fast for these largest examples, so that the system becomes too ill-conditioned before we reach optimality. For the specific example CANT-8-2-2-7, we set the initial $\epsMINRES=10^{-5}\sqrt{n}$, because this additional accuracy was required for convergence. Such non-default parameter choices are indicated in Tables~\ref{tab:CANT_large} and \ref{tab:BRIDGE_large}.

It needs to be said that even with such adjustments, the PBM algorithm did not solve all problems to the specified accuracy. For all BRIDGE-$m_x$-$m_y$-$m_z$-$\ell$ examples with $\ell=7$, it failed either close to or in the last iteration, after $\delta(u,\alpha) / \frac{1}{2}f^\top u$ had dropped below $\epsPBM=10^{-5}$ and $\epsNWT$ had been set to $10^{-4}$. The residual term $\resPBM$, defined in \eqref{eq:resIP}, did not go below $\epsNWT$ as required. This was because at a certain point, the approximate solutions of the reduced Newton system \eqref{eq:system} were no longer directions of descent, presumably due to numerical errors. In these cases, we accepted the, as it were, nearly optimal solutions at which the algorithm stalled. The iteration numbers we list in the table are those after which no further change in residual values is seen. Note that $\resPBM$ was well below $10^{-3}$ for all cases except BRIDGE-8-2-2-7, and the scaled duality gap of the accepted solution was always below $10^{-5}$.

It is evident from Tables~\ref{tab:CANT_large} and \ref{tab:BRIDGE_large}, that the PBM method is both more efficient and more robust than the IP method. In both cases, the use of a multigrid preconditioner for the MINRES solver achieves the desired result in that the number of MINRES iterations grows sublinearly with the size of the system, if at all. The CANT-$m_x$-$m_y$-$m_z$-$\ell$ example even displays a decrease in MINRES iterations with larger system size in some cases. However, this is probably not representative and a possible explanation involves the parameter $\epsMINRES$: since its initial value scales with the problem size, it might simply be chosen lower than necessary for the smaller problems.

\section{Conclusion}\label{sec:conclusion}
In this paper, we proposed a PBM method to solve the dual of the VTS formulation of the minimum compliance topology optimization problem. We compared it with the DOC method, one of the most popular methods for topology optimization, on the one hand, and with the IP method as an established method for general convex problems, on the other. The implementations of both the PBM and IP algorithms were tailored to the specific problem. All three methods used a multigrid preconditioned MINRES solver for the linear systems arising in each iteration.

In our numerical experiments, the PBM method clearly came out on top. It was around 20 times faster in terms of CPU time than the OC method when requiring the same degree of optimality. Even when using a very generous stopping criterion in the OC method---one that yields visibly sub-optimal results---PBM was still faster.

The IP method suffers from the characteristic ill-conditioning of the system matrix, which in some of our experiments prevented convergence altogether. Here, PBM proved to be much more robust, in addition to being considerably faster. Still, convergence was not guaranteed for all large-scale examples when sticking to the strictest stopping criterion. Judging by the symmetry and smoothness of the final design, the results were still satisfactory. Overall, the convergence behavior of the PBM method seems to be sensitive to changes in parameters such as stopping tolerances or scaling parameters. A thorough parameter study might further improve the algorithm.

We did not consider the DOC method for such large-scale problems, as its expected computation time simply disqualified it as a competitor. It is however possible that it would eventually converge even for those problems where PBM does not. Note that this would most likely take days or even weeks, as compared to the typical (successful) PBM run which took less than 12 hours. Since the DOC does not feature multipliers or barrier-/penalty-parameters tending to 0, it is not as susceptible to ill-conditioning as the PBM or IP method. This means that the advantage of DOC, when compared with PBM, could be reliability, albeit at the price of serious inefficiency.
\section*{Acknowledgements}\label{sec:Ack}
The authors would like to thank Michael Stingl for the use of the Matlab routines used to visualize the results in Section~\ref{subsec:m-scale}.

\bibliography{VTS_IP_PBM_MGM}
\end{document}